\documentclass[11pt]{amsart}

\usepackage{lineno}
\usepackage{hyperref}
\usepackage[top=4cm, bottom=4cm, left=3cm, right=3cm]{geometry}
\usepackage{amsmath} 
\usepackage{amsthm} 
\usepackage{mathrsfs}
\usepackage{amsfonts}
\usepackage{mathtools}
\usepackage{amssymb}
\usepackage{xcolor}
\usepackage{graphicx}

\usepackage{cite}
\usepackage{url}

\makeatletter
\newcommand{\pcite}[2][]{%
	(\kern-0.35em
	\@ifempty{#1}{\cite{#2}}{\cite[#1]{#2}}%
	)%
}
\makeatother

\usepackage[draft]{todonotes}   

\newtheorem{thm}{Theorem}[section]
\newtheorem{corollary}[thm]{Corollary}
\newtheorem{proposition}[thm]{Proposition}
\newtheorem{theorem}[thm]{Theorem}
\newtheorem{lemma}[thm]{Lemma}

\theoremstyle{definition}

\newtheorem{example}[thm]{Example}
\theoremstyle{remark}
\newtheorem{remark}[thm]{Remark}
\numberwithin{equation}{section}

\newcommand{\RR}{\mathbb{{R}}}
\newcommand{\ZZ}{\mathbb{{Z}}}
\newcommand{\NN}{\mathbb{{N}}}
\newcommand{\TT}{\mathbb{{T}}}
\newcommand{\CC}{\mathbb{{C}}}
\newcommand{\DD}{\mathbb{{D}}}

\newcommand{\linearOp}{\mathcal{L}}

\newcommand\Real{{\mathfrak R}{\mathfrak e}\,} %

\usepackage{hyperref}
\hypersetup{pdfpagemode={UseOutlines},	bookmarksopen, pdfstartview={FitH},	colorlinks,	linkcolor={red},	citecolor={blue},	urlcolor={blue}}

\makeatletter
\newcommand{\setword}[2]{%
	\phantomsection
	#1\def\@currentlabel{\unexpanded{#1}}\label{#2}%
}
\makeatother

\begin{document}
	\title[On similarity to contraction semigroups and tensor products, I]{On similarity to contraction semigroups\\ and tensor products, I}
	
	\author{J. Oliva-Maza}
	\address[J. Oliva-Maza]{Departamento de Matem\'aticas, Instituto Universitario de Matem\'aticas y Aplicaciones,
		Universidad de Zaragoza, 50009 Zaragoza, Spain}
	\email{joliva@unizar.es}
	\thanks{The first author has been mainly supported by XXXV Scholarships for Postdoctoral Studies, by Ram\'on Areces Foundation. He has also been partially supported by Project PID2022-137294NB-I00, DGI-FEDER, of the MCYTS and Project E48-23R, D.G. Arag\'on, Universidad de Zaragoza, Spain.
	}
	
	\author{Y. Tomilov}
	\address[Y. Tomilov]{Institute of Mathematics, Polish Academy of Sciences, {\'S}niadeckich 8, 00-956 Warsaw, Poland, and Faculty of Mathematics and Computer Science, Nicolas Copernicus University, Chopin Street 12/18, 87-100 Toru{\'n}, Poland}
	\email{ytomilov@impan.pl}
	
	\thanks{The second author was supported by the NCN Opus grant UMO-2023/49/B/ST1/01961. He was also partially supported by the NAWA/NSF grant BPN/NSF/2023/1/00001 and the NCN Weave-Unisono grant
		2024/06/Y/ST1/00044.}

	
	
	
\begin{abstract}
In the context of finite tensor products of Hilbert spaces, 
we prove that similarity of a tensor product of operator semigroups to a contraction semigroup 
is equivalent to the corresponding similarity for each factor, after an appropriate rescaling. 
A similar result holds with contractivity replaced by quasi-contractivity. 
This splitting phenomenon allows us to construct new and, in a sense, the strongest possible examples of 
$C_0$-semigroups not similar to contractions, thus completing an important chapter of the theory. 
We also address the discrete setting and relate it to our results.
\end{abstract}

	\maketitle

\tableofcontents

\section{Introduction}

\subsection{State of the art, results and their motivations}\label{simil0}

Let $H$ be a Hilbert space and $\linearOp(H)$ be the space of bounded linear operators on $H.$
The study of similarity  to a contraction on
$H,$ and more generally to a family of contractions on $H,$
is a classical and deep subject of operator theory,
with many connections to other domains of analysis.
It goes
back to classical results of Rota and  Sz.-Nagy, who noted, respectively,
that
if $T$ is a bounded operator on $H$ with spectral radius  strictly less
than $1$ or $T$ is compact and power bounded,
then $T$ is similar to a contraction on $H$.
Moreover, Sz.-Nagy proved that if $T$ is doubly power bounded, i.e.
$\sup_{n \in \mathbb Z}\|T^n\|<\infty,$ then $T$ is even similar to a unitary operator on $H.$

After the results of Sz.-Nagy and Rota
there was a hope that every power bounded
operator is similar to a contraction,
which was refuted by Foguel's example in \cite{foguel1964counterexample}.
Soon after Halmos reworked it  in \cite{halmos1964foguel}. Moreover,
since Foguel's operator appeared to be not polynomially bounded \pcite{lebow1968power-bounded}, he strengthened
the hypothesis and asked in \cite{halmos1970ten} whether every polynomially bounded
operator on $H$ is similar to a contraction.
This question attracted a considerable attention, giving rise to a separate line of research.

Without going into details of numerous attempts
to obtain a positive answer, we emphasize a paper related to the subject  of this paper.
Paulsen, Pearcy and Petrovi{\'c} proved in \cite{paulsen1995centered} that if $T\in \linearOp(H)$ is centered, that is, $ T^m T^{*m}$ and $T^{*n} T^n$  commute
for all $m, n \in \mathbb N,$
then the polynomial boundedness of $T$ implies
its similarity to a contraction.
The arguments in \cite{paulsen1995centered} relied in particular on the fact (see Theorem $2.3$ there)
that if
$T_1 \in \linearOp(H)$ is arbitrary and $T_2$ has norm and spectral radius equal to $1,$
then the tensor product $T_1\otimes T_2$ is similar to a contraction on $H\otimes H$ only if
$T_1$ is so, and the same holds for polynomial boundedness instead of similarity to a contraction.
In fact, for the similarity claim the norm constraint can be dropped:
if both $T_1$ and $T_2$ have spectral radius $1,$
then similarity of $T_1
\otimes T_2$ to a contraction
is equivalent to the same property satisfied separately by $T_1$ and $T_2.$
Although the
statement has not yet found significant applications so far, the methodology of \cite{paulsen1995centered} is of interest,
and it can play a role in the study of semigroups, as this paper shows.

After several milestone results, in particular, by Peller \cite{peller1982estimates}, Bourgain \cite{bourgain1986similarity},  Paulsen \cite{paulsen1984every},
and Aleksandrov and Peller \cite{aleksandrov1996hankel},
a seminal counterexample to Halmos' question was found by Pisier in \cite{pisier1997polynomially}, see also \cite {davidson1997polynomially}  and \cite{pisier2001similarity}
for supplementary results and a pertinent discussion.
Other related and notable examples were elaborated in \cite{kalton2002solution} and \cite{pisier2001multipliers}.

In the light of this counterexample and with aim of understanding the multivariate
setting, it is natural  to address the problem of \emph{joint similarity} to contractions
and to ask whether  for commuting operators $T_1, T_2 \in \linearOp(H)$
the property $\|S_1T_1S_1^{-1}\|\le 1$ and $\|S_2T_2S_2^{-1}\|\le 1$
for some invertible $S_1, S_2 \in \linearOp(H)$
implies  that
$\|ST_1S^{-1}\|\le 1$ and $\|ST_2S^{-1}\|\le 1$
for an invertible $S \in \linearOp(H).$
It was shown by Petrovi\'c in \cite{petrovic1997polynomially}
that one can find  a polynomially bounded $T\in \linearOp(H)$ resembling Pisier's example, such that
$(T\otimes I)(I \otimes T)=T \otimes T$ is not polynomially bounded,
thus the product of commuting polynomially bounded operators need not be polynomially bounded.
This example was substantially improved by Pisier in \cite{pisier1998joint}, who
obtained a negative answer to the joint similarity problem. His counterexample is a clever adaptation of Foguel's type examples
depending on  ideas from \cite{bozejko1987littlewood} and \cite{peller1982estimates}.

Despite the (counter-)examples by Pisier and Petrovi\'c,
a variety of positive results, stating  joint similarity to contractions for a finite family of commuting operators,
was obtained through the last two decades.
In particular, analogues of Rota's theorem for commuting families of operators can be found
in \cite{fong1978renorming} and \cite{shih1984analytic}. 
Among many other results addressing joint similarity to contractions, one  can mention criteria for  finite families of commuting power bounded matrices
\cite{clouatre2019joint},
families of Foguel-Hankel type operators with operator entries \cite{constantin2010joint}, \cite{constantin2009joint}, 
families of, in general, non-commuting operators subject to appropriate conditions \cite{popescu2011joint},  \cite{popescu2014similarity}, and families of 
commuting Ritt operators \cite{arrigoni2021square}, \cite{arrigoni2019functional}, 
  \cite{mohanty2019joint}.

The research on joint similarities to contractions on $H$ fits naturally into the study of $C_0$-semigroups
similar to semigroups of contractions on $H$, denoted by $\mathcal{SC}(H)$, in the sequel.
First, $\mathcal T=(T(t))_{t \ge 0} \in \mathcal{SC}(H)$
means that there exists a semigroup of contractions $\mathcal S=(S(t))_{t \ge 0}$ and an invertible operator $R \in \linearOp(H)$
such that  $RT(t)R^{-1}=S(t), t \ge 0,$  and so   $(T(t))_{t \ge 0}$ are jointly similar to $(S(t))_{t \ge 0}.$
Second, it is of interest to explore joint similarity of a family of commuting $C_0$-semigroups
to a family of contraction semigroups, when in the definition above $\mathcal T$ and $\mathcal S$ may vary but $R$ remains the same.
The present paper addresses both kinds of joint similarity and link them to each other.
As we will see, the continuous case is rather different from the discrete one.

The study of semigroups in $\mathcal{SC}(H)$ has developed in parallel with the discrete case. 
It was mainly motivated by applications to PDE theory via the Lumer-Phillips theorem. 
Following Foguel's example, Packel produced in \cite{packel1969semigroup} a bounded $C_0$-semigroup not similar to a semigroup of contractions. 
This example was improved by Benchimol in \cite{benchimol1977feedback}, who used a construction similar to Packel's, although the technical details were quite different. 
Benchimol's semigroup had an explicitly given generator, making it amenable to further modifications. 
In an important paper \cite{chernoff1976two}, Chernoff constructed a $C_0$-semigroup on $H$ which is not similar to a semigroup of contractions even after rescaling, i.e.  multiplying by $e^{-at}$ for any $a \in \mathbb{R}$, and thus provided an exponentially stable $C_0$-semigroup with this property. 
Note that for $T \in \linearOp(H)$, the similarity of $T^n$ to a contraction on $H$ for some $n \in \mathbb{N}$ implies that $T$ itself is similar to a contraction, and thus $(T^n)_{n \ge 0}$ is jointly similar to contractions on $H$. 
Already \cite{chernoff1976two} showed that an analogue of this property for a $C_0$-semigroup $\mathcal{T} = (T(t))_{t \ge 0}$ does not in general hold, providing additional insights into the joint similarity properties of $\mathcal{T}$. 
Perhaps the most striking example of a semigroup outside $\mathcal{SC}(H)$ was obtained by Le Merdy \cite{lemerdy2000bounded}, relying on ideas from \cite{baillon1991examples}. 
The example yields a $C_0$-semigroup $\mathcal{T} = (T(t))_{t \ge 0}$ on $H$ such that $\mathcal{T}$ is not similar to a contraction semigroup, and at the same time $\mathcal{T}$ is exponentially stable, sectorially bounded holomorphic of angle $\pi/2$, and $T(t)$ is compact for all $t > 0$. 
Thus, in view of Chernoff's and Le Merdy's examples, continuous analogues of Nagy's and Rota's criteria fail dramatically.

The only hope was to strengthen exponential stability by nilpotency with or even without
compactness, and this question was asked explicitly in \cite[p. 423]{vu1997almost} and in \cite[p. 202]{vu1998similarity}.
In this paper, we refute this conjecture
as well as its variants.
Our approach to the similarity problems rests on the notion of quasi-contractivity, which
is of independent interest. Recall that a semigroup $\mathcal T=(T(t))_{t\ge 0} \subset \mathcal L(H)$ is said to be quasi-contractive if there exists $a \in \mathbb R$ such that the rescaled semigroup
\[
e_{a}\mathcal T:=(e^{at}T(t))_{t \ge 0}
\]
is contractive.
 The class 
$\mathcal{SQC}(H)$ of semigroups on a Hilbert space 
 that are similar to a quasi-contraction semigroup is of fundamental importance in semigroup theory and in its applications to PDEs, and it appears in a variety of contexts. Although seemingly close to 
$\mathcal{SC}(H)$, the class 
$\mathcal{SQC}(H)$ is much larger and includes, for example, all 
$C_0$-groups as well as $C_0$-semigroups generated by bounded perturbations of normal operators.
 On the other hand, Chernoff's
example says that there exists a bounded $C_0$-semigroup, which is not similar even to a
quasi-contraction semigroup. Moreover, it is not difficult to show that Le Merdy’s semigroup
is also not similar to a quasi-contraction semigroup.

\subsection{Tensor products of semigroups and similarity problems}\label{simil}

The aim of this paper is two-fold.
First,
for an integer $m \ge 2$ and $k=1, \dots, m,$ 
given a Hilbert space $H_k$ and a nonzero semigroup $\mathcal T_k=(T_k(t))_{t \ge 0} \subset L(H_k)$ strongly continuous on $(0,\infty),$ 
we consider the tensor product semigroup $\otimes_{k=1}^m \mathcal T_k =(\otimes_{k=1}^m T_k(t))_{t \ge 0}$ defined on the tensor product $\otimes_{k=1}^{m}H_k$ of the $H_k$'s.
We elaborate an interesting splitting phenomenon 
saying that $\otimes_{k=1}^m \mathcal T_k$ is similar to a contraction (or quasi-contraction)  semigroup
if and only if each of the factors $\mathcal T_k$ is so, possibly after rescaling.
This property can be viewed as a continuous, multivariate counterpart of the
Paulsen-Pearcy-Petrovi\'c theorem stated in the preceding section.
Note that there are many statements in the literature on preservation of a certain property by tensor products.
However, the converse results are rare, more involved, and do not usually concern metric properties.

Setting for convenience $\NN_m:=\{1, \dots,  m\},$ we now give a precise formulation of our result.
\begin{theorem}\label{split_finite}
	Let $m \ge 2$ be an integer. For every $k \in \NN_m$, let $H_k$ be a Hilbert space 
	and let $\mathcal T_k = (T_k(t))_{t\geq0} \subset \mathcal L(H_k)$
	be a nonzero semigroup, strongly continuous in $(0,\infty)$. Then
	\begin{itemize}
		\item	[(i)]  $ \otimes_{k=1}^m \mathcal T_k \in \mathcal{SQC} \left( \otimes_{k=1}^m H_k \right)$ if and only if $\mathcal T_k \in \mathcal{SQC}(H_k)$
		for all $k \in \NN_m$.
		\item [(ii)] 
		$ \otimes_{k=1}^m \mathcal T_k \in \mathcal{SC} \left(\otimes_{k=1}^m H_k\right)$ if and only if there exist $(d_k)_{k=1}^m \subset \RR$ with $\sum_{k=1}^m d_k = 0$
		such that $e_{d_k}\mathcal T_k \in \mathcal{SC}(H_k)$
		for all $k \in \NN_m$.
	\end{itemize}
\end{theorem}

Using a variant of Rota's theorem for quasi-contraction semigroups, we relate the classes 
$\mathcal{SQC}(H)$ and 
$\mathcal{SC}(H),$ and, employing quasi-contractivity as a bridge, 
deduce Theorem~\ref{split_finite}(ii)  from Theorem~\ref{split_finite}(i).

Furthermore, we present our results for semigroups that are not necessarily strongly continuous at $0$. Although such semigroups are relatively uncommon in the context of Hilbert spaces, our proofs readily extend to this general setting, as they do not rely on the generator of the semigroup. 
Moreover,
it appeared that the splitting phenomenon is so strong, that it persists, under natural technical assumptions, even in the case of  infinite
tensor products. This much more demanding setting will be addressed in a separate paper 
\cite{oliva2025tensorbis}. In the studies of infinite tensor products of semigroups,
minimal regularity assumptions become crucial, especially in the situation of semigroups
defined on ``very large'', complete infinite tensor products
 $\otimes_{k=1}^\infty H_k$ of Hilbert spaces $H_k$
 which apart from degenerate situations cannot accommodate even unitary $C_0$-groups.

It appears that the tensor product structure in Theorem~\ref{split_finite} is indispensable.
Observe that a  tensor product $\otimes_{k=1}^m \mathcal T_k$  of semigroups $\mathcal T_k, 1 \le k \le m,$ fits to the more general situation of a product
of commuting semigroups
$$\mathcal T_1\otimes \ldots \otimes I, \ldots, I\otimes \ldots \otimes \mathcal T_m,
$$
on $\otimes_{k=1}^m H_k$. 
However, in the general case of commuting $C_0$-semigroups on $H$
the splitting phenomenon fails dramatically, demonstrating that our setup
is natural, and
cannot be substantially generalized.
More precisely, we show that neither $\mathcal {SC}(H)$ nor $\mathcal{SQC}(H)$ is invariant under products of commuting semigroups.
This yields, in particular, a semigroup analogue of Pisier's example from \cite{pisier1998joint}.

\begin{theorem}\label{commutingCounter_intro}
	There exist  Hilbert spaces $H$ and $K$ and pairs of commuting $C_0$-semigroups $(\mathcal T_1, \mathcal T_2)$ and $(\mathcal S_1, \mathcal S_2)$ on $H$ and $K$ respectively
	such that
	\begin{itemize}
		\item [(i)] $\mathcal T_1,$ $\mathcal T_2 \notin \mathcal{SQC}(H)$, while their product $\mathcal T_1 \mathcal T_2$ belongs $\mathcal{SC}(H);$
		\item [(ii)] $\mathcal S_1, \mathcal S_2  \in \mathcal{SC}(K),$ yet their  product
		$\mathcal S_1 \mathcal S_2$ lies outside $\mathcal {SQC}(K).$
	\end{itemize}
\end{theorem}

The proof of Theorem~\ref{commutingCounter_intro}(i) relies on Benchimol's example and Chernoff's direct sums trick from \cite{chernoff1976two}.
To prove Theorem~\ref{commutingCounter_intro}(ii) we employ Pisier's construction of two bounded operators that are similar to contractions but whose product is not even polynomially bounded. Additionally, we use an interpolation technique originating from Bhat-Skeide (see \cite{bhat2015pure}), with further elaboration on these ideas provided in \cite{dahya2024interpolation}.

\subsection{Tensor products of operators and similarity problems}\label{simil2}

We also prove discrete variants of the above splitting results. These are less demanding, as they avoid the significant difficulties that arise from dealing with $T(t)$ for small $t$, a problem that plays no role in the discrete setting. Adapting the proof of the Paulsen--Pearcy--Petrovi\'c theorem from \cite{paulsen1995centered}, which concerns the discrete setting, we obtain a discrete analogue of Theorem \ref{split_finite}~(ii). By characterizing semigroups $\mathcal{T} = (T(t))_{t \ge 0}$ in $\mathcal{SC}(H)$ via bounds for the similarity constants $\mathcal{C}(T(t))$ (Proposition~\ref{CrSimSetProp}), we then deduce Theorem~\ref{split_finite}~(ii)
from \cite{paulsen1995centered}, except in the quasi-nilpotent case. The Paulsen--Pearcy--Petrovi\'c theorem itself relies on Paulsen's deep completely boundedness criterion for similarity to a contraction, which we avoid entirely. In this approach, similarity to a quasi-contraction semigroup must be handled separately and still does not cover quasi-nilpotent semigroups, making it unsuitable for constructing examples of quasi-nilpotent $C_0$-semigroups outside $\mathcal{SC}(H)$ (see Subsection~\ref{simil3}). For these reasons, we argue directly, in order to address issues specific to continuous time.

Apart from \cite{paulsen1995centered}, tensor product considerations similar in spirit to our work can be found in \cite{gupta2019neumann} and \cite{kubrusly1997introduction}, which also address the discrete setting. In \cite[Proposition 8.9]{kubrusly1997introduction}, it was shown that a left shift $S_T$ on $\ell^2(\mathbb{N}, H)$ with operator weight $T \in \linearOp(H)$ is similar to a contraction if and only if $T$ is. Hence, as noted in \cite[p.~118]{kubrusly1997introduction}, choosing $T$ to be, for example, Foguel's operator, then $S_T$ is a power-bounded operator that satisfies $S_T^n \to 0$ strongly, but $S_T$ is not similar to a contraction. Clearly, $S_T = S \otimes T$.  
For $n$-tuples $\mathcal{A}$ and $\mathcal{B}$ in $\linearOp(H^n)$, it was proved in \cite[Theorem 1.2]{gupta2019neumann} that if $\mathcal{B}$ satisfies a matrix version of the (multivariable) von Neumann inequality and contains $1$ in its (Harte) spectrum, then $\mathcal{A} \otimes \mathcal{B}$ satisfies von Neumann's inequality if and only if the inequality holds for $\mathcal{A}$.  
The statements and arguments in these works are relatively straightforward, whereas our approach relies on more advanced methods and techniques and establishes more general results.

\subsection{Examples of semigroups not similar to contraction ones}\label{simil3}

To achieve the second aim of our work, we employ Theorem~\ref{split_finite} to produce new examples of semigroups not similar to contractions, thereby filling gaps in the theory. As mentioned
in Section~\ref{simil0}, several examples are now known in the literature demonstrating that continuous versions of basic criteria for the similarity of $T \in \mathcal{L}(H)$ to a contraction can fail dramatically. In particular, Le~Merdy's example of a $C_0$-semigroup $\mathcal{T} = (T(t))_{t \ge 0}$ outside $\mathcal{SC}(H)$ shows that neither a version of Nagy's condition nor a version of Rota’s condition holds for $\mathcal{T}$, even though $T(t)$ satisfies the corresponding discrete versions of these conditions for every $t>0$. This supports the conjecture that a certain uniformity in the similarity properties of $T(t)$ might imply $\mathcal{T} \in \mathcal{SC}(H)$, a claim confirmed by Proposition~\ref{CrSimSetProp} on similarity constants. However, such uniformity cannot be obtained by considering only $t \ge \tau$ for some fixed $\tau > 0$, since it is present in the aforementioned examples (see \cite[Section 7.3]{oliva2025similarity}). Several results emphasizing the importance of the behaviour of $(T(t))_{t\ge 0}$ near zero for similarity problems can be found in \cite[Section 4]{oliva2025similarity}.

Thus, the only remaining hope of preserving even traces of Nagy's and Rota's theorems was to impose very strong assumptions, such as the nilpotency or quasi-nilpotency of $\mathcal{T} = (T(t))_{t \ge 0}$, possibly combined with the compactness of $T(t)$ for $t > 0$ (semigroups satisfying the latter property are called immediately compact in the literature). In this setting, the behaviour of $T(t)$ for large $t$ is as good as possible, while severe restrictions are imposed on its behaviour near zero. Since all known examples of $\mathcal{T}$ outside $\mathcal{SC}(H)$ have generators with at least a countable spectrum, they can hardly be adapted directly to produce new examples. As a result, the existence of a nilpotent or even quasi-nilpotent semigroup not similar to a semigroup of contractions remained an open problem for some time (see \cite[p. 426]{vu1997almost} for an explicit question and further remarks). To overcome these and other difficulties, we employ tensor products of semigroups together with Theorem~\ref{split_finite}. This approach allows us to construct which seem to be the strongest possible examples of $\mathcal T$ lying outside $\mathcal{SC}(H)$, obtained as straightforward corollaries of Theorem~\ref{split_finite}.

We begin by constructing a nilpotent $C_0$-semigroup $\mathcal{T} = (T(t))_{t \ge 0}$ on $H$ with $T(t)$ compact for every $t > 0$, such that $\mathcal{T} \notin \mathcal{SC}(H)$. 
It is well known that nilpotent operators admit a canonical triangular representation (see, e.g., \cite{andruchow1988nilpotent}), which greatly simplifies their study. 
Since we are not aware of any analogous representation for continuous-parameter semigroups, we follow an alternative approach based on tensor products. 
Furthermore, removing the nilpotency assumption allows us to add a substantial amount of regularity and to produce a quasi-nilpotent, compact, and holomorphic $C_0$-semigroup of maximal angle $\pi/2$, which is again not similar to a semigroup of contractions. 
(Observe that a nilpotent semigroup cannot be holomorphic.)
Thus the following statement is true.
\begin{theorem}\label{Counter_intro}
There exist $C_0$-semigroups $\mathcal{T}$ and $\mathcal{S}$ on a Hilbert space $H$ such that neither belongs to $\mathcal{SQC}(H)$, with $\mathcal{T}$ nilpotent and immediately compact, and $\mathcal{S}$ quasi-nilpotent, immediately compact, and bounded holomorphic of angle $\pi/2$.
\end{theorem}
Theorem \ref{Counter_intro} reveals that no direct analogue of the discrete similarity conditions can guarantee similarity to a semigroup of contractions without additional assumptions on the behaviour of the semigroup near zero.

On the other hand, direct analogies do occur on the level of examples. Using the Bhat--Skeide interpolation technique together with Theorem~\ref{split_finite}, we construct a variety of bounded semigroups outside $\mathcal{SC}(H)$, thus obtaining continuous counterparts of known discrete results.  
The following is a sample of obtained results.
\begin{theorem}\label{MTTh1}
There exists a $C_0$-semigroup $\mathcal{T}=(T(t))_{t \ge 0}$ on a Hilbert space $H$ such that, for all $h \in H$, $h \neq 0$,
\begin{equation}
\lim_{t\to \infty} T(t)h = 0 \quad \text{weakly}, \qquad \inf_{t>0} \|T(t)h\|>0,
\end{equation}
and there is no $R \in \mathcal L(H)$ with zero kernel and dense range satisfying
\begin{equation}\label{inter}
T(t)R=RS(t),\qquad t \ge 0,
\end{equation}
for a contraction $C_0$-semigroup $(S(t))_{t\geq0}$ on $H$.
\end{theorem}
Thus, there exists a bounded $C_0$-semigroup $\mathcal{T} = (T(t))_{t \ge 0}$ on $H$ that cannot even be intertwined with a contraction semigroup.
Moreover, as we show in Section~\ref{MTSubsect}, $\mathcal{T}$ can be chosen to have a norm bound arbitrarily close to $1$, along with some additional desirable properties.

The above theorem and other results from Section~\ref{MTSubsect} provide continuous analogs to the results in \cite{kalton2002solution} and \cite{muller2007quasisimilarity}. In particular, in terms of \eqref{inter}, Theorem~\ref{MTTh1} yields a substantially stronger example than, e.g., examples in \cite{benchimol1977feedback}, \cite{kalton2002solution}, and \cite{packel1969semigroup}.

As an illustration of our discrete results, we present an explicit example of a power bounded operator $T$ on $H$ that is not similar to a contraction, yet satisfies $T^n \to 0$ and $T^{*n} \to 0$ strongly as $n \to \infty$. A similar example was originally obtained in \cite{eckstein1972contre-exemple} via a modification of Foguel's example. Another example, based on tensor products, shows the existence of a compact, injective, and quasi-nilpotent operator $K$ and a power bounded operator $A$ such that $K$ is in the bicommutant of $A$ and $A$ is not similar to a contraction. This provides a counterexample to the question posed in \cite[Problem 3]{vu1997almost}, now with the additional property that $K$ is quasi-nilpotent. A weaker counterexample could be derived from \cite{lemerdy2000bounded}.

It remains to note that the arguments in this paper are soft and rather non-technical relying
on simple auxiliary constructions. 

	\subsection{Notation} \label{NotationSection}
	
Here we fix notation that will be used throughout the paper.
All Hilbert spaces considered in this paper are assumed to be nonzero. Given two Hilbert spaces $H$ and $K$, $\linearOp(H,K)$ denotes the space of bounded operators from $H$ to $K$. If $H=K$, we write $\linearOp(H)$. $I_H$ denotes the identity operator on $H$, or simply $I$ when the choice of space is apparent, and we denote the identity semigroup $(I)_{t\geq0}$ by $\mathcal I$. 
Given a closed operator $A$ on $H$, $\sigma(A)$ and $\rho(A)$ stand for the spectrum  and  the resolvent set of $A,$ respectively. If $A$ is bounded, then $r(A)$ denotes the spectral radius of $A$.

Given a semigroup $\mathcal T = (T(t))_{t\geq0} \subset \linearOp(H)$,  recall that we write $\mathcal T \in \mathcal{SC}(H)$ if $\mathcal T$ is similar to a contraction semigroup, and $\mathcal T \in \mathcal{SQC}(H)$ if $\mathcal T$ is similar to a quasi-contraction semigroup. If $\mathcal T \in \mathcal{SC}(H)$, then $\mathcal C(\mathcal T)$ denotes the similarity constant of $\mathcal T$. Analogously, if $A\in \linearOp(H)$ is similar to a contraction, then $\mathcal C(A)$ denotes the similarity constant of $A.$ The exponential growth bound of $\mathcal T$ is denoted by $\omega_0(\mathcal T)$, and
for $d \in \RR$, $e_d \mathcal T$ stands for the semigroup $(e^{dt} T(t))_{t\geq0}$.

We use the symbol $\otimes$ to denote the tensor product $\otimes_{k=1}^m H_k$ of Hilbert spaces 
$(H_k)_{k=1}^m$, the tensor product $\otimes_{k=1}^m A_k$ of operators 
$(A_k)_{k=1}^m$, or the tensor product $\otimes_{k=1}^m \mathcal T_n$ of semigroups
 $(\mathcal T_k)_{k=1}^m$.

Banach limits on $\ell^\infty(\mathbb N)$ will be denoted by ${\rm LIM}.$
We will not distinguish between different Banach limits, 
and will write ${\rm LIM} \bigl[ x_k \bigr]$ for $(x_k)_{k=1}^{\infty} \in \ell^\infty(\mathbb N).$

With a slight abuse of notation, as is customary in the literature, we denote the norm of a bounded operator $A$ on a Hilbert space with a fixed norm $\|\cdot\|$ by $\|A\|$. This helps avoid overloaded notation when the space is considered with multiple norms.

We use $\simeq$ to stand for a unitary equivalence of either spaces or operators.
The  notation $\times_{i \in I} G_i$ and its relatives will denote the direct products of linear 
spaces $G_i, i \in I.$

Finally, we set $\mathbb N_m=\mathbb N \cap\{1, \dots, m\}$ for $m \in \mathbb N,$ $m\geq 2$, $\CC^+ := \{ z\in \CC \, : \, \Real z > 0\}$, and $\TT := \{ z \in \CC \, : \, |z| =1\}.$

\section{Preliminaries of semigroups and relevant tools}\label{Prelim}
		
	In this section, we introduce and discuss several basic notions pertaining to tensor products, semigroups and similarity to contractions. Although most of this material is standard, we present it here  
	in concise form, for clarity and ease of reference
	
		A family $\mathcal T = (T(t))_{t\geq0}$ of bounded operators on a Hilbert space $H$ is said to be a \textit{semigroup} (of bounded operators) if $T(0)=I$ and  $T(s+t) = T(s)T(t)$ for $s,t\geq0$.
		 Recall that a semigroup is strongly continuous in $(0,\infty)$ if and only if it is strongly measurable, that is, for each $h \in H$, the mapping from $[0,\infty)$ to $H$, given by $t\mapsto T(t)h$, is Bochner measurable, see for instance \cite[Theorem 10.2.3]{hille1957functional}. Also, if for each $h\in H$, the mappings $t \mapsto T(t)h$ are norm continuous in $[0,\infty)$, then $\mathcal T$ is said to be a $C_0$-semigroup. This last property is equivalent to $\lim_{t\to0^+} T(t)h = h$ for all $h \in H$.

	Given a semigroup $\mathcal T = (T(t))_{t\geq0}$ on a Hilbert space $H$ that is strongly continuous in $(0,\infty)$, we denote by $\omega_0(\mathcal T)$ its exponential growth bound, defined as
	$$\omega_0(\mathcal T) = \lim_{t\to \infty} \frac{\log \|T(t)\|}{t}  \in [-\infty, \infty).
	$$
	Note that $\omega_0(\mathcal T) = \frac{\log r(T(t))}{t}$ for every $t>0$, where $r(\cdot)$ denotes the spectral radius and $\log(0):=-\infty$. If moreover $\limsup_{t\to0} \|T(t)\| <\infty$ (this is always the case if $\mathcal T$ is a $C_0$-semigroup), then 
	$\mathcal T$ is locally bounded on $(0,\infty)$ and
	$$\omega_0(\mathcal T) = \inf\{\omega_0 \in \RR \, : \, \exists K\geq 1 \, \mbox{such that}\,  \|T(t)\| \leq K e^{\omega_0 t}, \, t \geq 0\}.
	$$ 
The next simple result on lower bounds for semigroup orbits will be crucial for the sequel.

\begin{lemma}\label{vanNerveenLemma}
\begin{itemize}
\item [(i)]	Let $\mathcal T = (T(t))_{t\geq0}\subset \linearOp(H)$ be a semigroup on a Hilbert space $H$, strongly continuous in $(0,\infty)$, such that $\limsup_{t\to0} \|T(t)\| < \infty$ and $\omega_0(\mathcal T) \geq 0$. Then
for every $t_0>0$ there exists a unit vector $x \in H$ satisfying 
 $\|T(t)x\| \geq 1/2$ for all $t \in [0,t_0].$
\item [(ii)]	Let $T \in \mathcal L(H)$ with $r(T)\ge 1.$ Then for every $n_0\in\mathbb N$ there exists a unit vector $x\in H$
such that $\|T^n x\|\ge 1/2$ for all $n \in  \NN_{n_0}.$
	\end{itemize}
\end{lemma}
The proof of a statement more general than (i) can be found in \cite[Lemma 3.1.7]{van1996asymptotic},
where the argument given for $C_0$-semigroups holds verbatim for semigroups bounded near $0$.
In fact, the proof of (i) relies on
(ii), which can be found, again in a more general version, in \cite[Theorem V.37.8]{muller2007spectral}. While the proof of (ii) is rather direct, its variant for (i) seems to be more involved,
and we stated both (i) and (ii)
for ease of reference.
	
	We say that  bounded operators $A$ and $B$ on a Hilbert space $H$ are \textit{similar} if there exists an invertible $R \in \linearOp(H)$ such that $A  = R B R^{-1}$. Analogously,  semigroups of bounded linear operators
	$ 
\mathcal T = (T(t))_{t\geq 0}$ and $\mathcal S = (S(t))_{t\geq 0}$ 
	are said to be \textit{similar} if $T(t) = R S(t) R^{-1}$ for all $t \geq 0$ and for some invertible $R \in \linearOp(H)$. If in addition $\sup_{t\geq 0} \|S(t)\| \leq 1,$ then we say that $\mathcal T$ is similar to a contraction semigroup.
			We denote by $\mathcal{SC}(H)$ the class of semigroups similar to contraction ones. 
	Recall that $A \in \mathcal{L}(H)$ is similar to a contraction if and only if  there exists an inner product $\langle \cdot, \cdot \rangle_{\mathscr H} : H \times H \to \CC$ such that the induced norm $\|\cdot\|_{\mathscr H}$ is equivalent to $\|\cdot\|$ and $\|Ax\|_{\mathscr H}\le \|x\|_{\mathscr H}$
	for all $x \in H.$ Similarly,
	$\mathcal T \in \mathcal{SC}(H)$ if and only if  there exists 
	a Hilbertian norm equivalent to the original one 
	such that 
	$$\|T(t)x\|_{\mathscr H} \le \|x\|_{\mathscr H}, \qquad t\geq 0.
	$$ 
	See e.g. \cite[pp. 235-236]{benchimol1977feedback} for a discussion of these well-known properties.
	Given a semigroup $\mathcal T = (T(t))_{t\geq0}$ in $\mathcal{SC}(H)$,
	we define the similarity constant $\mathcal C( \mathcal T)$ of $\mathcal T$ as
	\begin{equation*}
		\begin{aligned}
			\mathcal C(\mathcal T) :=
			  \inf \Big\{ \|R\| \|R^{-1}\| \, : \, R \in \linearOp(H) \mbox{ invertible with } \|R T(t) R^{-1}\| \leq 1 \mbox{ for all } t \geq 0\Big\},
		\end{aligned}
	\end{equation*}
	and observe that
	\begin{equation*}
		\begin{aligned}
			\mathcal C(\mathcal T) :=& \inf \Big\{ \mathcal C \, :  \, \mbox{ there exists a Hilbertian norm } \|\cdot\|_{\mathscr H} \mbox{ on } H \mbox{ such that }
			\\ & \qquad \qquad  \|T(t)\|_{\mathscr H} \leq 1, \, t\geq0, \mbox{ and } \|h\| \leq \|h\|_{\mathscr H} \leq \mathcal C \|h\| \mbox{ for all } h \in H \Big\},
					\end{aligned}
	\end{equation*}
%
	If $\mathcal T \notin \mathcal{SC}(H)$, we set $\mathcal C(\mathcal T) := \infty$. If $\mathcal T \in \mathcal{SC}(H)$, then there exists an equivalent Hilbertian norm $\|\cdot\|_{\mathscr H}$ on $H$ such that $\|h\| \leq \|h\|_{\mathscr H} \leq \mathcal C(\mathcal T) \|h\|$ for all $h \in H$ and $\|T(t)\|_{\mathscr H} \leq 1$ for all $t\geq0,$ i.e. $\mathcal C(\mathcal T)$ is attained. This was proven in \cite[Proposition 2.4]{holbrook1977distortion} for the similarity constant $\mathcal C(A)$ of an operator $A$ similar to a contraction, defined similarly as
	\begin{align*}
		\mathcal C(A) :=& \inf\{\mathcal C \, :  \, \mbox{ there exists a Hilbertian norm } \|\cdot\|_{\mathscr H} \mbox{ on } H \mbox{ such that }
		\\ \qquad & \|A\|_{\mathscr H} \leq 1,  \mbox{ and } \|h\| \leq \|h\|_{\mathscr H} \leq \mathcal C \|h\| \mbox{ for all } h \in H\},
	\end{align*}
	and the proof works with trivial modifications also in the semigroup setting,
	see \cite[p. 230]{holbrook1977distortion}, and cf. 
\cite[pp. 235-236]{benchimol1977feedback}.

The following generalization of contraction semigroups will be crucial for our purposes. A semigroup $\mathcal T = (T(t))_{t\geq0}$ on $H$ is said to be a \textit{quasi-contraction semigroup} if there exists $a \in \RR$ such that
	$$ \|T(t)\| \leq e^{a t}, \qquad t \geq 0,
	$$
	or equivalently, if $e_{-a}\mathcal T = (e^{-at} T(t))_{t\geq0}$ is a contraction semigroup. Then a semigroup $\mathcal T$ is said to be \textit{similar to a quasi-contraction semigroup} if there exists a quasi-contraction semigroup $\mathcal S = (S(t))_{t\geq0}$ on $H$ similar to $\mathcal T$, and we denote the class of quasi-contraction semigroups on a Hilbert space $H$ by $\mathcal{SQC}(H)$. It is clear that $\mathcal T \in \mathcal{SQC}(H)$ if and only if there exists $a\in \RR$ such that $e_{a} \mathcal T \in \mathcal{SC}(H)$.
	
In the study of similarity properties, as is often the case in this area, we will rely on the notion of Banach limits. This notion is indispensable in various averaging constructions and will be 
primarily used to average a sequence of auxiliary norms and construct an appropriate renorming.
Recall that for any abelian semigroup $S$ with the discrete topology
there exists a positive functional ${\rm LIM}$ on the space of bounded functions 
$\ell^\infty(S)$ satisfying ${\rm LIM}(1)=1$ and invariant with respect to left shifts
on $\ell^\infty(S).$
We will only need them for $S=\mathbb Z_+$, although $S=\mathbb R_+$ is also possible and leads to equivalent constructions.

\section{Tensor products of Hilbert spaces}\label{InfiniteTensorProdDefSect}
While tensor products of Hilbert spaces and their operators became a basic and well-understood
concept in functional analysis, 
 it is not easy to find their comprehensive treatment in the literature.
Yet  \cite{murray1936rings} and \cite{von1939infinite}  addressing finite and infinite tensor products respectively remain unsurpassable,
and since the case of finite tensor products trivially embeds into this setting, the paper  \cite{von1939infinite}
 still serves as a basic reference. A good elaboration of these works can be found in \cite{russmann2020tensor}.

Let $m \in \mathbb N, m \ge 2,$ be fixed, and let 
$\odot$ stand for algebraic tensor product of vector spaces. 
Given
a family of Hilbert spaces $(H_k)_{k=1}^m$
their tensor product $\otimes_{k=1}^m H_k$ 
is defined as a completion of 
$\odot_{k=1}^{m} H_k$ 
 under the canonical inner product:
\begin{equation}\label{tens_def}
\langle \odot_{k=1}^m f_k, \odot_{k=1}^m g_k \rangle = \prod_{k=1}^m \langle f_k,  g_k\rangle, 
\end{equation}
defined initially on elementary tensors $\odot_{k=1}^m f_k,   \odot_{k=1}^{m} g_k$ from $\odot_{k=1}^m H_k,$
and extended by linearity and density to $\otimes_{k=1}^m H_k$. 
In fact, one may omit (in a sense) the concept of algebraic tensor product and proceed from elementary tensors identified with multi-sesquilinear forms  $(f_k)_{k=1}^m \to \prod_{k=1}^m \langle f_k,  g_k\rangle,$
for all $
  (f_k)_{k=1}^m \in \times_{k=1}^m H_k.$
	Passing then to their linear span and completing it under the inner product
	defined as in \eqref{tens_def} and
extended to the span by linearity, we arrive at $\otimes_{k=1}^m H_k.$
These and other constructions lead to the same space $\otimes_{k=1}^m H_k$
since the properties \eqref{tens_def} and the density of the linear span of elementary 
tensors define $\otimes_{k=1}^m H_k$ uniquely, see e.g. \cite[Section 1]{russmann2020tensor}.
The properties of finite tensor products of Hilbert spaces are well-understood and can be found in many sources, although usually with poor content and very few details.
For comparatively complete (and complementing) treatments see  \cite[Section 3 and Section 4.6]{hackbush2019tensor},
\cite[p. 125-147]{kadison1983fundamentals}, \cite{murray1936rings}, 
\cite[Section 2.4]{reed1980functional} and \cite[Section 1]{russmann2020tensor}.
The approach in \cite[p. 125-147]{kadison1983fundamentals} is comprehensive, although somewhat non-standard. It reduces to the standard one as remarked in \cite[Section 1.3]{russmann2020tensor}.

It is often useful to know that if $(e_{j,k})_{j \in J_k}$ is an orthonormal basis of $H_k$ for $k \in \mathbb N_m,$ then 
$\{e_{j_1, 1}\otimes ... \otimes e_{j_m,m}: j_k \in J_k, \, k \in \mathbb N_m\}$ is an orthonormal basis 
of $\otimes_{k=1}^m H_k,$  see e.g. \cite[Lemma 2.2.1]{murray1936rings} or \cite[Proposition 2]{reed1980functional}.
Among other natural properties enjoyed by finite tensor products
we recall the \emph{associativity law} stating, in particular, that 
$\otimes_{k=1}^m H_k$ is (canonically) unitarily isomorphic
 with $H_j \otimes \left( \otimes_{k \neq j, k \in \mathbb N_m} H_k \right)$
for every $j \in \mathbb N_m,$ cf. \cite[Proposition 2.6.5]{kadison1983fundamentals}, \cite[Theorem VII]{von1939infinite}.

The next simple property of equivalent norms on $\otimes_{k=1}^m H_k$ 
will be key for proving similarity to contractions of various operator tensor products.
\begin{lemma}\label{moreEquivalentTensorRemark}
For each $k \in \mathbb N_m$, let $H_k$ be a Hilbert space and let $\| \cdot \|_{\mathscr{H}_k}$ be an equivalent Hilbertian norm on $H_k$ satisfying
	$$\|h\|_{H_k} \leq \|h\|_{\mathscr{H}_k} \leq \mathcal C_k \|h\|_{H_k}, \qquad h \in H_k,
	$$
		for some $\mathcal C_k >0.$ Then
	\begin{align}\label{similarNormTensorEq}
		\| \cdot\|_{\otimes_{k=1}^m H_k} \leq \| \cdot\|_{ \otimes_{k=1}^m \mathscr{H}_k} \leq \left(\prod_{k=1}^m \mathcal C_k\right) \| \cdot\|_{\otimes_{k=1}^m H_k}.
	\end{align}
	\end{lemma}
\begin{proof}
It suffices to note that for every $j \in \mathbb N_m$, using the associativity of finite tensor
 products of Hilbert spaces, we have
	\[
	\otimes_{k=1}^m H_k \simeq H_j \otimes \left( \otimes_{k \neq j, k \in \mathbb N_m} H_k \right) 
	\simeq \ell^2(I_j, H_j),
	\]
	where $\simeq$ stand for (canonical) unitary isomorphisms, and
$I_j$ indexes an orthonormal basis of  $\otimes_{k \neq j, k \in \mathbb N_m} H_k.$ 
Thus, 
	$$\| \cdot\|_{\otimes_{k=1}^m H_k} \leq \| \cdot\|_{\mathscr{H}_j  \otimes \left( \otimes_{k\neq j, k \in \mathbb N_m} \mathscr H_k\right)} \leq \mathcal C_j \| \cdot\|_{\otimes_{k=1}^m H_k}, \qquad j \in \mathbb N_m.
	$$
	An iteration of this argument yields 
		\begin{align*}\label{similarNormTensorEq}
		\| \cdot\|_{\otimes_{k=1}^m H_k} \leq 
		\| \cdot\|_{ \otimes_{k=1}^m \mathscr{H}_k} \leq \left(\prod_{k=1}^m \mathcal C_k\right)
		\| \cdot\|_{\otimes_{k=1}^m H_k}.
	\end{align*}	
\end{proof}

\section{Tensor products of operators and semigroups}\label{tensor_op}
 
The formula \eqref{tens_def} suggests a way to define operators on $\otimes_{k=1}^m H_k,$
and it appears to be the right route to proceed. 
Given  $(T_k)_{k=1}^m \subset \times_{k=1}^m\linearOp(H_k),$ 
there exists   the unique bounded operator $\otimes_{k=1}^m T_k$ on $\otimes_{k=1}^m H_k$ satisfying
\begin{equation}\label{def_fin}
(\otimes_{k=1}^m T_k) (\otimes_{k=1}^m f_k) = \otimes_{k=1}^m T_k f_k, \qquad \otimes_{k=1}^m f_k\in \otimes_{k=1}^m H_k,
\end{equation}
called the tensor product of $(T_k)_{k=1}^m.$
For more detailed information on finite Hilbert space tensor products of operators
one may consult \cite[p.144-147]{kadison1983fundamentals}, \cite[Chapter VIII.10]{reed1980functional}
or \cite[Chapters 3,4]{hackbush2019tensor}.
It is important to emphasize that by associativity of finite tensor products of Hilbert spaces,
a similar associative law holds for products of operators. In fact,
the associativity of finite tensor products of Hilbert spaces implies that the operator
$\otimes_{k=1}^m T_k$ is unitarily equivalent to
 $T_j \otimes \left( \otimes_{k \neq j, k \in \mathbb N_m} T_k \right)$ for every $j \in \mathbb N_m,$
via a canonical isomorphism.

Tensor products $\otimes_{k=1}^m T_k$
enjoy a number of special properties that distinguish them from general products $\prod_{k=1}^{m} T_k$ of commuting operators $T_k.$.
It is well-known that 
\begin{equation}\label{prod_norm}
\| \otimes_{k=1}^m T_k\| = \prod_{k=1}^m \|T_k\|,
\end{equation}
 see, for example, 
\cite[p. 146]{kadison1983fundamentals}
or \cite[p. 299-300]{reed1980functional}. 
Moreover, the spectra of the tensor products 
can be expressed in terms of their factors:
\begin{equation}\label{prod_spect_sp}
\sigma \left( \otimes_{k=1}^m T_k \right) = \prod_{k=1}^m \sigma(T_k).
\end{equation}
For $m=2$ with $H_1=H_2$, a proof of this well-known result can be found in \cite{brown1966spectra}, while the case of different $H_1$ and $H_2$ is treated in \cite[Theorem~XIII.34]{reed1978analysis}. The general case then follows from the associativity of the tensor products.

 For the context of Banach spaces and further generalizations one may consult e.g. \cite[Theorem 4.3]{ichinose1970spectra}.
The property \eqref{prod_spect_sp} implies an analogous relation for the spectral radii $r(T_k)$ of $T_k:$
\begin{equation}\label{prod_spect}
r\left( \otimes_{k=1}^m T_k \right) = \prod_{k=1}^m r(T_k),
\end{equation} 
which will be of value for the sequel. (Alternatively, to obtain \eqref{prod_spect}, one may use \eqref{prod_norm} along with Gelfand's 
spectarl radius formula.)
Note that \eqref{prod_spect} generalizes to the setting of infinite tensor products of operators.
The proof of this fact, given in \cite{oliva2025tensorbis},
depends on more involved notions and techniques.

Similarly to the discrete setting, for $m \in \mathbb N,$
given Hilbert spaces $(H_k)_{k=1}^m,$
and  one-parameter semigroups
$(\mathcal T_k)_{k=1}^m \subset \times_{k=1}^m \linearOp(H_k),$
their tensor product
 $\otimes_{k=1}^m \mathcal T_k =(\otimes_{k=1}^m T_k(t))_{t \ge 0}$ is defined 
on elementary tensors $\otimes_{k=1}^m f_k$ from $\otimes_{k=1}^m H_k$ by
\begin{equation}\label{def_fin_sem}
(\otimes_{k=1}^m T_k(t)) (\otimes_{k=1}^m f_k) = \otimes_{k=1}^m T_k(t) f_k 
\end{equation}
for every $t\ge 0$. 
Clearly, $\otimes_{k=1}^m \mathcal T_k$ is a one-parameter semigroup in
 $\mathcal L\left(\otimes_{k=1}^m H_k\right),$
and if $\mathcal T_k$ is strongly continuous on $(0,\infty)$ for every  $k \in \mathbb N_m,$ 
then $\otimes_{k=1}^m \mathcal T_k$ is
strongly continuous on $(0,\infty)$ as well,
see e.g. \cite[A-I, 3.7]{arendt1986one},  \cite[Sections 4,5]{gill1978infinite} for more details.

The following result may be well-known, 
but we could not find a reference for it in the literature.
\begin{lemma}\label{C0TensorLemma}
	For $m \ge 2$, let $(H_k)_{k=1}^m$ be a family of Hilbert spaces
	and let $(\mathcal T_k)_{k=1}^m   \subset \times_{k=1}^m \linearOp(H_k)$ be
	a family of semigroups 
	strongly continuous in $(0,\infty)$.
	Then $\otimes_{k=1}^{m} \mathcal T_k$ is a $C_0$-semigroup on 
	$\otimes_{k=1}^{m} H_k$ if and only if $\mathcal T_k$ is a $C_0$-semigroup on 
	$H_k$ for every $k \in \mathbb N_m.$
\end{lemma}
\begin{proof}
	
	It is direct to prove that if  $\mathcal T_k$ is a $C_0$-semigroup on $H_k$ for every $k \in \mathbb N_m$
	then $\otimes_{k=1}^{m} \mathcal T_k$ is a $C_0$-semigroup on $\otimes_{k=1}^{m} H_k$.
	(Note that $\otimes_{k=1}^{m} \mathcal T_k$ is a product of commuting strongly continuous
	semigroups.)
	
	To prove the opposite statement, assume first that $m=2.$
		Suppose that $\mathcal T_1 \otimes \mathcal T_2$ is a $C_0$-semigroup on $H_1\otimes H_2$, thus $\limsup_{t \to 0} \|T_1(t) \otimes T_2(t)\| < \infty$. Since $t\mapsto \|T_2(t)\|$ is a submultiplicative measurable on $(0,\infty),$ we infer that $\liminf_{t\to 0} \|T_2(t)\| \geq 1$, see for instance \cite[Theorem 7.4.3]{hille1957functional}. Consequently, we obtain that
	$$\limsup_{t\to 0} \|T_1(t)\| = \limsup_{t\to 0} \frac{ \|T_1(t) \otimes T_2(t)\|}{\|T_2(t)\|} < \infty,
	$$
	and thus $\|T_1(t)\|\le M e^{\omega t}$ for some $M\ge 1$, $\omega \in \RR$ and all $t \ge 0.$
	Then following the proof of \cite[Corollary 2.2]{arendt2001approximation}
	(see also \cite[Proposition 2.1]{arendt2001approximation}),  we set
	$$N:= \cap_{t>0} \ker T_1(t), \qquad \text{and}\qquad K:= \{h \in H_1\, :\, \lim_{t\to0} T_1(t) h = h\},
	$$
	and infer that  $N$ and $K$ are $\mathcal T_1$-invariant subspaces of $H_1$ such that 
	$K \oplus N = H_1$, where $\oplus$ denotes a topological direct sum, and the restriction of $\mathcal T_1$ to $K$ is a $C_0$-semigroup.
	It is clear that $\lim_{t\to 0} T_1(t) h_1 \otimes T_2(t)h_2 = 0$ for every $h_1 \in N$ and $h_2 \in H_2$. By the strong continuity of $\mathcal T_1 \otimes \mathcal T_2$ at zero, we must have $N = \{0\}$, 
	that is, $H_1 = K$ and consequently $\mathcal T_1$ is a $C_0$-semigroup. By symmetry,  
	$\mathcal T_2$ is also a $C_0$-semigroup by what we already proved.	Taking in account the associativity of tensor products, the general case is a direct implication
	of the case $m=2$ considered above.
\end{proof}

\begin{remark}
	The claim given by Lemma \ref{C0TensorLemma} also holds if the assumption on the strong continuity of each $\mathcal T_k$ on $(0,\infty)$ is replaced by the assumption that least $m-1$ of $(\mathcal T_k)_{k=1}^m$  are strongly continuous on nontrivial subsets of $(0,\infty)$. Moreover, this new assumption is, in general, optimal. Otherwise,
	one may consider semigroups $\mathcal T_1=
	(e^{i\gamma(t)}I)_{t \ge 0},$ 
	$\mathcal T_2=({e^{-i\gamma (t)}}I)_{t \ge 0},$ and $\mathcal T_k=\mathcal I, k\ge 3,$ if $m \ge 3,$
	where $\gamma: \mathbb R_+ \to \mathbb R_+$ is additive and non-measurable. (Such $\gamma$ exists by Hamel bases.)
	
	We omit the proof of this modified version of Lemma \ref{C0TensorLemma} as the arguments involved differ substantially from those used in the original setting with strong continuity on $(0, \infty)$ and this claim has no direct application to the main content of the paper. 
	\end{remark}

Note that 
if for $m \in \mathbb N,$ we are given
a family of Hilbert spaces $(H_k)_{k=1}^m$  and a family of locally bounded semigroups
$(\mathcal T_k)_{k=1}^m \subset \times_{k=1}^m\mathcal L (H_k),$ 
then 
\begin{equation}\label{exptype}
\omega_0 (\otimes_{k=1}^m \mathcal T_k)=\sum_{k=1}^m \omega_0 (\mathcal T_k).
\end{equation}
This follows directly from \eqref{prod_spect} applied to $T_k(t)$ for fixed $t>0.$
The relation \eqref{exptype} will be used repeatedly in the sequel.

\section{Similarity to contractions and finite tensor products of semigroups}\label{ContinuousSection}

In this section, we establish one of the main results of the paper,  Theorem~\ref{split_finite}, which essentially states that the tensor product of a finite number of semigroups is in $\mathcal{SC}\left(\otimes_{k=1}^m H_k\right)$ (in $\mathcal{SQC} \left(\otimes_{k=1}^m H_k\right)$) if and only if, up to a suitable rescaling, each of the factors is in $\mathcal{SC}(H_k)$ (in $\mathcal{SQC}(H_k)$). To this aim
we first prove Theorem \ref{split_finite}(i), providing the statement for quasi-contraction semigroups and being of independent interest.
 Then using 
a Rota-type theorem (Proposition~\ref{contractiveStableProp}), we extend this result to contraction semigroups, thus completing the proof.

\begin{proof}[Proof of Theorem~\ref{split_finite}(i)]
	
	For each $k \in \mathbb N_m$, let $H_k$ be a Hilbert space and let $\mathcal T_k = (T_k(t))_{t\geq0} \subset \linearOp(H_k)$ be a nonzero semigroup strongly continuous in $(0,\infty)$. We will show  that $\otimes_{k=1}^m \mathcal T_k \in \mathcal{SQC}\left(\otimes_{k=1}^m H_k\right)$ is equivalent to 
	$\mathcal T_k \in \mathcal{SQC}(H_k)$ for every  $k \in \mathbb N_m$.
	
The proof of the ``if'' part is direct. If $\mathcal T_k \in \mathcal{SQC}(H_k)$ for 
each $k \in \mathbb N_m$ then,
	in view of Lemma~\ref{moreEquivalentTensorRemark}, we may assume without loss of generality that $\mathcal T_k$, $k \in \mathbb N_m,$ are quasi-contraction semigroups themselves. Then, recalling \eqref{prod_norm},
		we conclude that $ \otimes_{k=1}^m \mathcal T_k$ is a quasi-contraction semigroup on 
	$ \otimes_{k=1}^m H_k$, and the ``if'' part of the statement holds.

	We now turn to the ``only if'' implication, which we prove by induction on $m$.
	 Let first $m=2.$ Let $H_1$ and $H_2$ be Hilbert spaces, and $\mathcal T_1 = (T_1(t))_{t\geq0} \subset \linearOp(H_1)$ and $\mathcal T_2 = (T_2(t))_{t\geq0} \subset\linearOp(H_2)$ be nonzero semigroups  strongly continuous  on $(0,\infty)$ and satisfying 
		$\mathcal T_1  \otimes \mathcal  T_2 \in \mathcal{SQC}(H_1 \otimes H_2)$. 
	Then clearly $\limsup_{t\to0} \|T_1(t) \otimes T_2(t)\| < \infty.$
	Since the mapping $(0,\infty)\ni t\to\|T_2(t)\|$ is measurable and submultiplicative, using \cite[Theorem 7.4.3]{hille1957functional} we have
	$\liminf_{t\to0} \|T_1(t)\| \geq 1.$ 
	Hence
	$$\limsup_{t\to 0} \|T_2(t)\| = \limsup_{t\to 0} \frac{\|T_1(t)\otimes T_2(t)\|}{\|T_1(t)\|} < \infty.
	$$ 
	From here, by the submultiplicativity of the map $[0,\infty)\ni t\mapsto \|T_2(t)\|,$
	it follows that $\sup_{t\in[0,b]} \|T_2(t)\|< \infty$ for each $b>0$.
	 Similarly, $\limsup_{t\to0} \|T_1(t)\| < \infty$ and $\sup_{t\in[0,b]} \|T_1(t)\|< \infty$ for every $b>0$.
	
	Let now $\|\cdot\|_{\rm{eq}}$ be an equivalent Hilbertian norm on $H_1  \otimes H_2$ 
	such that
	$\mathcal T_1  \otimes \mathcal T_2$ is quasi-contractive in $(H_1  \otimes H_2, \|\cdot\|_{\rm{eq}}),$ 
	and denoting by $\|\cdot\|$ the original norm on $H_1\otimes H_2$
	we have $\|\cdot\| \leq \|\cdot\|_{\rm{eq}} \leq \mathcal C \|\cdot\|$ for some $\mathcal C\geq 1.$
	Letting $\langle \cdot , \cdot \rangle_{{\rm eq}}$ be the corresponding inner product on $H_1 \otimes H_2,$ fix $b \in (0,\infty)$ and $h\in H_1$ such that $\|h\|=1$ and $T_1(b)h \neq 0$, and introduce the inner product
	$\langle \cdot, \cdot \rangle_{\mathscr{H}_2}$ on $H_2$ by
	\begin{align*}
		\langle f, g \rangle_{\mathscr{H}_2} :=  
		\int_0^b \langle T_1(t) h \otimes f, T_1(t) h \otimes g \rangle_{\rm{eq}} \, dt, \qquad f, g\in H_2.
	\end{align*}
	Note that the  integral above is well defined since the mapping $t\mapsto T_1(t)  \otimes I$ 
	from $(0,\infty)$ to $\mathcal L(H_1  \otimes H_2)$ is strongly continuous and bounded near $0$.
	Moreover, if 
	\begin{equation}\label{def_d}
	D:= \int_0^b \|T_1(t)h\|^2 \,dt,
	\end{equation}
	then $D \in (0,\infty)$ since $T_1(b)h \neq 0$ and $\sup_{t\in(0,b)} \|T_1(t)\|<\infty$. 
	Hence, 
	\begin{align*}
		\|f\|_{\mathscr{H}_2}^2 &= \int_0^b \|T_1(t) h \otimes f\|_{\rm{eq}}^2 \, dt \leq \mathcal C^2 \int_0^b \|T_1(t) h \otimes f\|^2 \, dt
		= D \mathcal C^2  \|f\|^2, \qquad f \in H_2,
	\end{align*}
	and similarly
	\begin{align*}
		\|f\|_{\mathscr{H}_2}^2 &\geq D \|f\|^2, \qquad f \in H_2,
	\end{align*}
	so that $\|\cdot\|_{\mathscr H_2}$ is an equivalent norm on $H_2$. 
	Fix $a\in \mathbb R$ such that $e_{-a} \mathcal T_1  \otimes \mathcal T_2$ is a contraction semigroup on $(H_1 \otimes H_2, \|\cdot\|_{\rm{eq}})$
	and assume, without loss of generality, that $a \ge 0.$ Setting $M:=\sup_{t \in (0,b)}\|T_2(t)\|,$
	observe that for all $s\in (0,b)$ and $f \in H_2$ with
	$\|f\|_{\mathscr H_2}=1,$
	\begin{equation}\label{AuxQuasiEq2}
		\begin{aligned}
		\|T_2(s) f\|_{\mathscr{H}_2}^2 &=  
				\int_0 ^s \|T_1(t)h \otimes T_2(s)f\|_{\rm{eq}}^2 \, dt+
							\int_s ^b \|T_1(t)h \otimes T_2(s)f\|_{\rm{eq}}^2 \, dt \\
							&\le \int_0^{s} e^{2at} \|h \otimes T_2(s-t)f\|_{\rm{eq}}^2 \, dt 
							+\int_0^{b-s} \|T_1(s+t)h \otimes T_2(s)f\|_{\rm{eq}}^2 \, dt \\
							&\leq \mathcal C^2 \left(\sup_{t\in(0,b)} \|T_2(t)\|^2 \right)
							\int_0^{s} e^{2at} \|h \otimes f\|^2 \, dt+
		\int_0^{b-s} e^{2as} \|T_1(t) h \otimes f\|_{\rm{eq}}^2 \, dt\\
			 &\leq D^{-1} se^{2as} \mathcal C^2 M^2   + e^{2as}
			\leq e^{\left(2a+\frac{\mathcal C^2 M^2}{D} \right)s}.
		\end{aligned}
\end{equation}
For general $s\ge0$, write $s=nb+r$ with $n\in\mathbb N$ and $r\in[0,b)$.
By the semigroup property,
$\|T_2(s)\|_{\mathscr H_2}\le \|T_2(b)\|_{\mathscr H_2}^n\,\|T_2(r)\|_{\mathscr H_2}$,
and the exponential estimate in  \eqref{AuxQuasiEq2} extends to 
all $s\ge0$.
	Hence $\mathcal T_2$ is a quasi-contraction semigroup on $(H_2,\|\cdot\|_{\mathscr H_2}),$
	and therefore $\mathcal T_2 \in \mathcal{SQC}(H_2)$. A similar reasoning shows that $\mathcal T_1 \in \mathcal{SQC}(H_1),$ 
	so the statement holds for $m=2$.
	
	To run the induction, we let now $m>2$ and assume that the statement of the theorem is true for all $k \in \NN$ such that $k< m$. Let $\mathcal T_k$, $k=1, \ldots, m,$ be such that $\otimes_{k=1}^m \mathcal T_k \in \mathcal{SQC}\left(\otimes_{k=1}^m H_k\right)$. By the associative law for tensor products,  
	$ \otimes_{k=1}^m \mathcal T_k \simeq \left( \otimes_{k=1}^{m-1} \mathcal T_k \right)  \otimes \mathcal T_m$ canonically. From here, by the induction assumption for $m=2$, it follows that $ \otimes_{k=1}^{m-1} \mathcal T_k \in \mathcal{SQC} \left(\otimes_{k=1}^{m-1} H_k\right)$ and $\mathcal T_m \in \mathcal{SQC}(H_m)$.  Another application of the induction assumption, this time for $k = m-1,$ yields that $\mathcal T_k \in \mathcal{SQC}(H_k)$ for $k=1,\ldots, m-1$. Hence the statement is true for $k=m$,  and the proof is completed.
\end{proof}

\begin{remark}\label{forinfinite}
Let $\mathcal T_1$ and $\mathcal T_2$ satisfy the assumptions of Theorem \ref{split_finite}(i).
If
 $\mathcal T_1 \otimes \mathcal T_2 \in \mathcal{SC}(H_1 \otimes H_2)$
and $\mathcal C$ and  $D$ are defined as in the proof of the theorem, then using \eqref{AuxQuasiEq2},
we conclude that
\begin{equation}\label{remarkk}
\|T_2(t)\|_{\mathscr{H}_2}\le e^{ \frac{\mathcal C^2 M^2}{2D} t}, \qquad t \ge 0,
\end{equation}
where $\mathcal C=\mathcal C(\mathcal T_1 \otimes \mathcal T_2)$ and $M=\sup_{t \in (0,b)} \|T_2(t)\|.$
Clearly, a similar estimate holds for $\mathcal T_1$. 
The estimate \eqref{remarkk} will play a crucial role in our studies of
infinite tensor products of semigroups in \cite{oliva2025tensorbis}.
\end{remark}

We proceed with Rota's theorem counterpart for semigroups in $\mathcal{SQC}(H)$. Apart from its auxiliary role in the proof of Theorem \ref{split_finite}, it is of independent interest and, importantly, requires no regularity assumptions on the semigroup.
The latter is of value in the study of semigroups on (complete) infinite tensor products 
where infinite tensor products of semigroups do not inherit any regularity of its factors. For alternative approach 
to this and similar results see \cite{oliva2025similarity}.

\begin{proposition}\label{contractiveStableProp}
	Let $H$ be a Hilbert space and let $\mathcal T = (T(t))_{t\geq0} \subset\linearOp(H)$ be a semigroup strongly continuous in $(0,\infty)$ such that $\mathcal T \in \mathcal{SQC}(H)$. Then, for each $a > \omega_0(\mathcal T)$, there exists an equivalent Hilbertian norm $\|\cdot\|_{\rm{eq}}$ on $H$ satisfying
	$$\|T(t)\|_{\rm{eq}} \leq e^{at}, \qquad t \geq 0.
	$$
\end{proposition}
\begin{proof}
	Replacing $\mathcal T$ by $e_{-d}\mathcal T = (e^{-dt} T(t))_{t\geq0},$  $a$ by $a-d$ for a suitable $d\geq 0$, and passing to an equivalent norm, we may assume that $\mathcal T$ is a contraction semigroup. 
	If $\omega_0(\mathcal T) = 0$, then 
	the statement is trivial since necessarily $a\geq 0$. 
	If $\omega_0(\mathcal T)<0,$ then fix $a \in (\omega_0(\mathcal T),0)$ and $b \in (\omega_0(\mathcal T),a)$, and let $K\geq 1$ be such that $\|T(t)\| \leq Ke^{bt}$ for $t \geq 0$. 
	
	For each $h\in H$, define 
	$$G_h(t):=\|T(t)h\|^2, \qquad t\geq0.
	$$
	One has that $G_h$ is non-increasing and $\lim_{t\to \infty} G_h(t)=0$, so $G_h$ 
	is of bounded variation on $[0,\infty).$ Then define the new norm $\|\cdot\|_{\rm{eq}}$ on $H$ as
	$$\| h \|_{\rm{eq}}^2 := G_h(0) - 2a \int_0^\infty e^{-2at} G_h(t) \, dt, \qquad h,g\in H.
	$$
	where the integral is convergent in view of $a>\omega_0(\mathcal T)$. 
		Moreover, for all $h \in H$ we have  $\|h\|_{\rm{eq}} \geq \sqrt{G_h(0)} = \|h\|$
		and
	$$\|h\|_{\rm{eq}}^2 \leq \|h\|^2 -2a  K^2 \|h\|^2 \int_0^\infty e^{2(b-a)t}\,dt = \left(1 - K^2\frac{a}{a-b}\right) \|h\|^2, \qquad h \in H.
	$$
	Hence, $\| \cdot \|_{\rm{eq}}$ is an equivalent norm on $H$. In addition, 
	it is readily seen that $\|\cdot\|_{\rm{eq}}$ satisfies the parallelogram law, hence it is a Hilbertian norm. 
	
	Now, integration by parts yields
	\begin{align*}
		\|h\|_{\rm{eq}}^2 &= G_h(0) -2a \int_0^\infty e^{-2at} G_h(t) \, dt = - \int_0^\infty e^{-2at} dG_h(t), \qquad h \in H.
	\end{align*}
	Therefore, taking into account that $G_h$ is nonincreasing, we have
	\begin{align*}
		\|T(s)h\|_{\rm{eq}}^2 &= - \int_0^\infty e^{-2at} dG_{T(s)h}  = - e^{2as} \int_s^\infty e^{-2at} dG_h(t) 
		\\&\leq - e^{2as} \int_0^\infty e^{-2at} dG_h(t) = e^{2as} \|h\|_{\rm{eq}}^2, \qquad s\geq 0, \, h \in H.
	\end{align*}
	 Thus,  $\|T(s)\|_{\rm{eq}} \leq e^{as}$ for all $s \geq 0,$
	and the statement follows.
	\end{proof}

Now, we are ready to prove one of the main results of this paper.

\begin{proof}[Proof of Theorem \ref{split_finite}(ii)]
	
	As in the proof of Theorem \ref{split_finite}(i), for every $k \in \mathbb N_m$, let $H_k$ be a Hilbert space and let $\mathcal T_k = (T_k(t))_{t\geq0} \subset \linearOp(H_k)$ be a nonzero semigroup strongly continuous in $(0,\infty)$. 
	We show that $\otimes_{k=1}^m \mathcal T_k \in \mathcal{SC}\left(\otimes_{k=1}^m H_k\right)$ if and only if there exist $d_k \in \RR$, $k=1,\ldots,m$ such that $\sum_{k=1}^m d_k=0$ and $e_{d_k} \mathcal T_k \in \mathcal{SC}(H_k),$ $k \in \mathbb N_m$.
	
	For the ``if'' part of the statement, using Lemma ~\ref{moreEquivalentTensorRemark}, we may assume that $e_{d_k} \mathcal T_k$ is a contraction semigroup on $H_k$ for every  $k \in \mathbb N_m.$
	Then
	\begin{align*}
		\left\| \otimes_{k=1}^m T_k(t)\right\| &= 	\left\| \otimes_{k=1}^m e^{d_k t} T_k(t)\right\| = \prod_{k=1}^m  \left\|e^{d_k t} T_k(t) \right\| \leq 1, \qquad t \geq 0,
	\end{align*} 
	and $ \otimes_{k=1}^m \mathcal T_k$ is a contraction semigroup on $ \otimes_{k=1}^m H_k$. \medskip
	
	Next, we show the ``only if'' part by induction in $m$. Let $m=2$, and assume that $\mathcal T_1  \otimes \mathcal T_2 \in \mathcal{SC}(H_1\otimes H_2)$. We then show that there exists $d\in \RR$ such that $e_{-d} \mathcal T_1 \in \mathcal{SC}(H_1)$ and $e_d \mathcal T_2 \in \mathcal{SC}(H_2)$.

Using \eqref{exptype},
	we have $ \omega_0(\mathcal T_1) + \omega_0(\mathcal T_2) = \omega_0(\mathcal T_1  \otimes \mathcal T_2) \leq 0$. If $\omega_0(\mathcal T_1) + \omega_0(\mathcal T_2) < 0$, then choose $d \in (\omega_0(\mathcal T_1), -\omega_0(\mathcal T_2))$ so that $\omega_0(e_{-d}\mathcal T_1) = \omega_0(\mathcal T_1) - d< 0$ and $\omega_0(e_d \mathcal T_2) = d + \omega_0(\mathcal T_2) < 0$. Since $e_{-d} \mathcal T_1 \in \mathcal{SQC}(H_1)$ and $e_d \mathcal
	T_2 \in \mathcal{SQC}(H_2)$ by Theorem~\ref{split_finite}(i), it follows from Proposition~\ref{contractiveStableProp} that $e_{-d} \mathcal T_1 \in \mathcal{SC}(H_1)$ and $e_d \mathcal T_2 \in \mathcal{SC}(H_2)$, as required.
	
	If $\omega_0(\mathcal T_1) + \omega_0(\mathcal T_2) = 0$,  then 
	replacing $\mathcal T_1$ and $\mathcal T_2$ with $e_{-d} \mathcal T_1$ and $e_d \mathcal T_2,$ respectively,  we may assume $\omega_0(\mathcal T_1) = \omega_0(\mathcal T_2) = 0$ where $d: = \omega_0(\mathcal T_1) = -\omega_0(\mathcal T_2)$. Since $\mathcal T_1 \otimes \mathcal T_2 \in \mathcal{SC}(H_1\otimes H_2)$ 
	we have $\sup_{t\geq0} \|T_1(t) \otimes T_2(t)\| < \infty.$
	Thus, taking into account that $\|T_1(t)\| \ge r(T_1(t)) = 1$ for all $t\ge 0$, we conclude that
	$$
	\sup_{t\geq0} \|T_2(t)\| = \sup_{t\geq0} \frac{\|T_1(t) \otimes T_2(t)\|}{\|T_1(t)\|} 
	\leq \sup_{t\geq0} \|T_1(t) \otimes T_2(t)\| < \infty.
	$$
	Similarly, we obtain that $\sup_{t\geq0} \|T_1(t)\| <\infty,$  and
	let
	\begin{equation}\label{bounded_sem}
	M:= \sup_{t\geq0} \left\{  \max \left\{\| T_1(t)\|, \|T_2(t)\|\right\} \right\}.
	\end{equation}
	By assumption, there exists an inner product $\langle \cdot,\cdot \rangle_{{\rm eq}}$ on $H_1\otimes H_2$    such that the corresponding norm 	$\| \cdot \|_{\rm{eq}}$ 	satisfies
	$$
	\|\cdot\| \leq \|\cdot\|_{\rm{eq}} \leq \mathcal C (\mathcal T_1\otimes \mathcal T_2) \|\cdot\|
	$$
	and $\mathcal T_1 \otimes \mathcal T_2$ is contractive
	on $H_1\otimes H_2$ with respect to $\|\cdot\|_{{\rm eq}}.$
		 From Lemma~\ref{vanNerveenLemma} it follows that there exists $(h_n)_{n=1}^\infty \subset H_1$ satisfying $\|h_n\| = 1$ and $\|T_1(t) h_n\| \geq 1/2$ for all $t\in [0,n],$ $n \in \NN$.
	Now, for each $n \in \NN$, define the inner product $\langle \cdot, \cdot \rangle_n$ on $H_2$ by
	\begin{align*}
		\langle f, g \rangle_n := \frac{1}{n} \int_0^n \langle T_1(t) h_n \otimes f, T_1(t)h_n \otimes g \rangle_{\rm{eq}} \, dt, \qquad f,g \in H_2.
	\end{align*}
	Note that the integral above is well-defined since the semigroup  $\mathcal T_1 \otimes \mathcal I$ is uniformly bounded and strongly continuous in $(0,\infty)$.  
	By the choice of $(h_n)_{n=1}^\infty$, we have
	$$
	D := \sup_{n\in \NN} \frac{1}{n} \int_0^n \|T_1(t) h_n\|^2 \, dt < \infty
	\qquad \mbox{and} \qquad d:= \inf_{n\in \NN} \frac{1}{n} \int_0^n \|T_1(t) h_n\|^2 \, dt > 0.
	$$
 	Moreover, denoting $\mathcal C:=\mathcal C (\mathcal T_1\otimes \mathcal T_2),$ 
	\begin{align}\label{upper_est}
		\|f\|_n^2 = \frac{1}{n} \int_0^n \|T_1(t) h_n \otimes f\|_{\rm{eq}}^2 \, dt \leq  \frac{\mathcal C^2}{n} \int_0^n \|T_1(t) h_n \otimes f\|^2 \, dt
				\leq  \mathcal C^2 D \|f\|^2, 
	\end{align}
	for all $f \in H_2$ and $n \in \NN.$
	Similarly, we conclude that
	\begin{align}\label{lower_est}
		\|f\|_n^2 &\geq   {d} \|f\|^2, \qquad f \in H_2, \, n \in \NN.
	\end{align}
	Next, taking into account \eqref{bounded_sem} and recalling that 
	$\| T_1(t) \otimes T_2(t) \|_{\rm{eq}}\le 1$ for all 	$t \ge0,$
		we observe that
	\begin{equation}\label{contr_est}
		\begin{aligned}
			\|T_2(s) f \|_n^2
	 		=&\frac{1}{n}\int_0^s \|T_1(t)h_n \otimes T_2(s)f\|_{\rm{eq}}^2 \, dt+\frac{1}{n}\int_s^n \|T_1(t)h_n \otimes T_2(s)f\|_{\rm{eq}}^2 \, dt\\
	\le &\frac{\mathcal C^2}{n} \int_0^s \left\|  
			h_n \otimes T_2(s-t) f \right\|^2\, dt + \frac{1}{n} \int_0^{n-s} \left\| T_1(t) h_n \otimes  f \right\|_{\rm{eq}}^2\, dt  \\
		\le&	\frac{\mathcal C^2}{n} \int_0^s \|T_2 (t)f\|^2\, dt + \|f\|_{n}^2
					\le \left(1  + \frac{\mathcal C^2 M^2 s}{dn} 
		\right) \|f\|_n^2
		\end{aligned}
		\end{equation}
		for all $n \in \NN$, $s \in (0,1)$ and $f \in H_2$.
	
	Finally, define the new inner product on $H_2$ by
	$$
	 \langle f, g \rangle_{\mathscr{H}_2} := \operatorname{LIM} \bigl[\langle f, g \rangle_n \bigr], \qquad f,g
	\in H_2,
	$$
	where $\operatorname{LIM}$ denotes a Banach limit on $\ell_{\infty}(\mathbb N).$ Combining \eqref{upper_est} and \eqref{lower_est}, we infer  that the corresponding norm $\| \cdot \|_{\mathscr{H}_2}$ is Hilbertian and equivalent to the original norm on $H_2.$
	Moreover, in view of \eqref{contr_est}, 
	\begin{align*}
		\|T_2(s) f\|_{\mathscr{H}_2} ^2 
		\leq  \operatorname{LIM}
		\left[ \left(1 +  \frac{\mathcal C^2 M^2 s}{d n}\right)   
		\|f\|_n^2 \right] = \|f\|^2_{\mathscr{H}_2}, 
	\end{align*}
	for all $s \in (0,1)$ and $f \in H_2$. Therefore, $\mathcal T_2 \in \mathcal{SC}(H_2)$, and by  an analogous argument,  $\mathcal T_1 \in \mathcal{SC}(H_1).$

We then proceed by induction similarly to the proof of Theorem \ref{split_finite}(i).	
Let $m>2$ and assume that the statement  is true for all $k \in \NN$ such that $k< m$. Let 
$\mathcal T_k$, $n=1, \ldots, m$, be such that $\otimes_{k=1}^m \mathcal T_k \in \mathcal{SC} \left(\otimes_{k=1}^m H_k\right)$. By the associativity of  tensor products, one has
 $\otimes_{k=1}^m \mathcal T_k \simeq \left( \otimes_{k=1}^{m-1} \mathcal T_k \right)  \otimes \mathcal T_m$ 
canonically.
Hence, by the induction assumption for
 $k=2$, there exists $d_m \in \RR$ such that $e_{-d_m}  \otimes_{k=1}^{m-1} \mathcal T_k \in \mathcal{SC}\left(\otimes_{k=1}^{m-1} H_k\right)$ and $e_{d_m} \mathcal T_m \in \mathcal{SC}(H_m)$. 
Moreover, the induction assumption for $k = m-1$ implies  that there is $(d_k)_{k=1}^{m-1} \subset  \RR$ satisfying $\sum_{k=1}^{m-1} d_k = - d_m$ and $e_{d_k} \mathcal T_k \in \mathcal{SC}(H_k)$ for every $k \in \mathbb N_{m-1}.$
Hence the statement holds for $k=m$ as well, and the proof is thus finished.

\end{proof}

	In the proof of Theorem~\ref{split_finite}(ii), if $\omega_0(\mathcal T_1 \otimes \mathcal T_2)  =0$, then
	for each $n \in \mathbb N$ one may define the equivalent inner product 
	$ \langle\cdot,\cdot\rangle_n'$ on $H_2$ as
	$$\langle f, g \rangle _n' := \frac{1}{D_n} \int_0^n
	\langle T_1(t) h_n \otimes f, T_1(t)h_n \otimes g\rangle_{\rm{eq}} \, dt, \qquad D_n:=\int_0^n 
	\| T_1(t)h_n \|^2\, dt,
	$$
		for all $f\in H_2$ and $n \in \NN.$ Then the corresponding  norm on $H_2$ satisfies
	$$\|f\| \leq \|f\|_n'\leq \mathcal C (\mathcal T_1\otimes \mathcal T_2) \|f\|, \qquad f \in H_2, \, n \in \NN,$$
	hence $\mathcal C (\mathcal T_2) \leq \mathcal C(\mathcal T_1\otimes \mathcal T_2)$.
	The argument used in the induction step from the proof above
	allows to extend the ``only if'' part of Theorem \ref{split_finite}(ii)  to the following
	form.
\begin{theorem}\label{splitbis}
	For each $k \in \mathbb N_m$, let $H_k$ be a Hilbert space and let $\mathcal T_k = (T_k(t))_{t\geq0} \subset \linearOp(H_k)$ be a nonzero semigroup strongly continuous in $(0,\infty)$ such that $\otimes_{k=1}^m \mathcal T_k \in \mathcal{SC}\left(\otimes_{k=1}^m H_k\right)$ and $\omega_0\left(\otimes_{k=1}^m \mathcal T_k\right) = 0$. Then there is $(d_k)_{k=1}^m \subset \RR$ with $\sum_{k=1}^m d_k = 0$ such that 
	$$
	\max_{k \in \mathbb N_m}\mathcal C (e_{d_k} \mathcal T_k) \leq \mathcal C \left(\otimes_{k=1}^m \mathcal T_k\right).
	$$ 
\end{theorem}	

	We do not know if the assumption $\omega_0\left(\otimes_{k=1}^m \mathcal T_k\right) = 0$ is redundant in the statement above.

To finish this section, we employ Theorem~\ref{split_finite}(ii) to characterize semigroups $\mathcal T$ such that $\mathcal T \otimes \mathcal S$ is similar to a contraction semigroup if and only if $\mathcal S$ is similar to a contraction semigroup. This characterization seems to be of value in applications.

	Let $H$ be a Hilbert space, and let $\mathcal T$ be a semigroup on $H$ strongly continuous in $(0,\infty)$. We say that $\mathcal T$ \textit{tensorially preserves similarity to contraction semigroups} if, for every Hilbert space $K$ and every semigroup $\mathcal S$ on $K$ strongly continuous in $(0,\infty)$, $\mathcal T  \otimes \mathcal S \in \mathcal{SC}(H \otimes K)$ is equivalent to $\mathcal S \in \mathcal{SC}(K)$.

\begin{proposition}\label{TensorPreserving}
	Let $H$ be a Hilbert space, and let $\mathcal T = (T(t))_{t\geq0} \subset \linearOp(H)$ be a semigroup strongly continuous in $(0,\infty)$. Then $\mathcal T$ tensorially preserves similarity to contraction semigroups if and only if $\mathcal T \in \mathcal{SC}(H)$ and $\omega_0(\mathcal T)=0$.
\end{proposition}
\begin{proof}
	Assume first that $\mathcal T$ tensorially preserves similarity to contraction semigroups, 
	and let $\mathcal I$ be the trivial identity semigroup on $\CC$. 
	As $\mathcal I \in\mathcal{SC}(\CC)$,  
	Using the canonical identification $H \simeq H \otimes \CC$ given by $h \mapsto h \otimes 1$, 
	one obtains $\mathcal{SC}(H  \otimes \CC)  \simeq  \mathcal{SC}(H)$ and 
	$\mathcal T  \otimes \mathcal I \simeq \mathcal T,$ so that $\mathcal T \in \mathcal{SC}(H)$ and, in particular, $\omega_0(\mathcal T) \leq 0$. Assume $\omega_0(\mathcal T) < 0$ and let $\mathcal S = e_d \mathcal I  \subset \linearOp(\CC)$, where $d \in \left(0, -\omega_0(\mathcal T)\right)$. Clearly, $\mathcal S \notin \mathcal{SC}(\CC),$ and  it 
	then follows that $e_d \mathcal T \simeq \mathcal T \otimes \mathcal S \notin \mathcal{SC}(H)$. However, as $\mathcal T \in \mathcal{SC}(H)$, one has $e_d \mathcal T \in \mathcal{SQC}(H)$ and $\omega_0(e_d \mathcal T) = \omega_0(\mathcal T) + d < 0$. Hence, Proposition~\ref{contractiveStableProp} implies that $e_d \mathcal T \in \mathcal{SC}(H)$, arriving at a contradiction. Thus we conclude that $\omega_0(\mathcal T) = 0$, and the 
	``only if'' part of the statement holds. \medskip
	
	To prove the ``if'' part, 
	we may assume that $\mathcal T$ is a contraction semigroup on $H$, see Remark~\ref{moreEquivalentTensorRemark}. As $\omega_0(\mathcal T) = 0$, we then have  $\|T(t)\| = 1$ for all $t\geq 0$. 
		Choose a semigroup $\mathcal S = (S(t))_{t\geq 0}$ on a Hilbert space $K,$  strongly continuous in $(0,\infty)$. If $\mathcal S \in \mathcal{SC}(K)$, then it is easy to see that $\mathcal T  \otimes \mathcal S \in \mathcal{SC}(H\otimes K).$ (This also follows from Theorem~\ref{split_finite}(ii) with $d_1 = d_2 =0$). 
	
	Now assume that $\mathcal T  \otimes \mathcal S \in \mathcal{SC}(H \otimes K)$. An application of Theorem~\ref{split_finite}(ii) shows that there exists $a\in \RR$ such that $e_a \mathcal S \in \mathcal{SC}(K)$ and $e_{-a} \mathcal T \in \mathcal{SC}(H)$. Since $\omega_0(e_{-a}\mathcal T) = -a + \omega_0(\mathcal T) =-a$, one obtains $a\geq 0$. Therefore, $\mathcal S \in \mathcal{SC}(K)$, as required.
	\end{proof}

To illustrate the result above, we mention the next immediate corollary.

\begin{corollary}\label{tensorShifts}
	Let $\mathcal T = (T(t))_{t\geq 0}$ be either  left or right shift semigroup  on one of the spaces
	$L^2(\Delta),$ where $\Delta$ is $[0,\infty), (-\infty, 0]$ or $\RR$. Then $\mathcal T$ tensorially preserves similarity to contraction semigroups. 
\end{corollary}
In particular, if $\mathcal R$ is the right shift on $L^2(\Delta)$, then using the identifications
 $L^2(\Delta)\otimes H\simeq L^2(\Delta, H)$ 
and $\mathcal T\otimes \mathcal R \simeq \mathcal T \mathcal R,$
we infer that the operator-weighted right shift  $\mathcal T \mathcal R$ belongs to $\mathcal{SC}(L^2(\Delta, H))$ if and only if 
the operator weight $\mathcal T$ belongs to 
$\mathcal{SC}(H).$
For an instance of utility of this fact see
\cite[Proposition 4.11]{oliva2025similarity}, where it arises in the study of  similarity properties of  Hilbert space semigroups, and is treated separately.
 
\section{When similarity to contractions fail}\label{counterexamplesSect}

Similarity to contraction semigroup may fail in various ways,
and in this section, using tensor products, we show that the failure could be rather dramatic.

\subsection{Commuting operators}\label{CommutingCounterSubsection}

We start by showing that neither of the implications in Theorem~\ref{split_finite}(ii) extends to the more general setting of commuting semigroups. 

First, we prove the failure of the ``only if'' claim, as stated in Theorem~\ref{commutingCounter_intro}(i). 
To provide the example required by Theorem~\ref{commutingCounter_intro}(i), we use as a building block a 
uniformly bounded $C_0$-semigroup $\mathcal S \in \mathcal{SQC}(H)\setminus \mathcal{SC}(H).$ 
To this aim, we 
 recall its construction
 given in \cite[Subsection 7.1]{oliva2025similarity} 
and relied on a modification of the argument in \cite{packel1969semigroup}. See also 
\cite{foguel1964counterexample} for an idea, which was elaborated further in \cite{packel1969semigroup}. Moreover,
variants of the construction, developed in \cite{oliva2025similarity}, allowed us to equip $\mathcal S$ with additional properties of independent interest

Let $J$ be either $\ZZ$, $\ZZ_+ = \{n \in \ZZ \, : \, n \geq 0\}$ or $\ZZ_- = \{n \in \ZZ \, : \, n \leq 0\}$, and fix $t>0$. If there exists $n \in J$ such that $3^n < t \leq 3^{n+1}$, then define $n_0(t) = n$. If there exists no such $n \in J$, then either $J = \ZZ_-$ and $t> 3$, and then we set $n_0(t)=0$; or $J = \ZZ_+$ and $t \leq 1$, and then we formally set $n_0(t) = -\infty$ so that $3^{n_0(t)}=0$. For each $n \in J$ such that $n > n_0(t)$, set $I_n(t) := (3^n-t, 3^n]$ and let $I_{n_0}(t) = [0, 2 \cdot  3^{n_0} -t]$ for $n_0= n_0(t)$. Thus, the operator $V_J(t) \in \linearOp(L^2(\RR_+))$ given by
\begin{align*}
	(V_J(t)f)(x): = \begin{cases}
		f  \left(2 \cdot 3^n - x - t \right), &\quad x \in I_n(t), \, n\geq n_0(t), 
		\\ 0, &\quad \mbox{otherwise},
	\end{cases}
	\qquad \qquad x \in (0,\infty), \,f \in L^2(\RR_+),
\end{align*}
is well-defined. We also set $V_J(0) = 0$ and define the family $\mathcal T_J = (T_J(t))_{t\geq0} \subset \linearOp (L^2(\mathbb R_+)\oplus L^2(\mathbb R_+))$ as
\begin{equation*}
	T_J(t): = \begin{pmatrix}
		L(t) & V_J(t)
		\\0 & R(t)
	\end{pmatrix},
	\qquad t\geq0,
\end{equation*}
where   $\mathcal L = (L(t))_{t\geq0}$  and $\mathcal R = (R(t))_{t\geq0}$ are the left shift semigroup 
and right shift semigroup, respectively, on $L^2(\RR_+)$. It was shown in \cite[Theorem 7.1]{oliva2025similarity} that $\mathcal T_J$ is a uniformly bounded $C_0$-semigroup on $L^2_{\oplus}(\mathbb R_+) := L^2(\mathbb R_+)\oplus L^2(\mathbb R_+)$ satisfying the following.
\begin{enumerate}
	\item [(i)] For every $t>0,$ the operator $T_J(t)$ is similar to a contraction if and only if $J= \ZZ_-$.
	\item [(ii)] $\mathcal T_J \in \mathcal{SQC}(L^2_{\oplus}(\mathbb R_+))$ if and only if $J = \mathbb Z_+$.
\end{enumerate}
Thus $\mathcal S = \mathcal T_{\ZZ_+} \in \mathcal{SQC}(L^2_{\oplus}(\mathbb R_+))\setminus \mathcal{SC}(L^2_{\oplus}(\mathbb R_+)).$

As an alternative to $\mathcal T_{\ZZ_-},$ one may consider 
a semigroup constructed by Benchimol in \cite[Section 4.3]{benchimol1977feedback}, following the same ideas as  \cite{packel1969semigroup}, and being motivated by 
\cite{foguel1964counterexample}.  Note, however, that the example in \cite[Section 4.3]{benchimol1977feedback} makes the generator explicit, which is often helpful in applications. 
More precisely, let $P$ be the orthogonal projection from $L^2(\RR_+)$ onto $L^2(\Delta)$, where
$$
\Delta := \bigcup_{n \in \NN} [3^n - 1, 3^n].
$$
and let $\mathcal D$ be the generator of the right shift semigroup $\mathcal R$ on $L^2(\RR_+)$, so that
its adjoint $\mathcal D^\ast$ generates the left shift semigroup $\mathcal L$ on the same space. Fix $\varepsilon > 0$, and set
$$A := \begin{pmatrix}
		\mathcal D^\ast & \varepsilon P \\
		0 & \mathcal D.
	\end{pmatrix}
$$
Then $A$ generates a quasi-contraction $C_0$-semigroup $\mathcal S$, since $A$ is a bounded perturbation of the generator $\begin{pmatrix} \mathcal D^\ast & 0 \\ 0 & \mathcal D \end{pmatrix}$ of the contraction $C_0$-semigroup $\mathcal L \oplus \mathcal R$. From an argument similar to Foguel's one, it follows that
  $\sup_{t \geq 0} \|S(t)\| \leq 1 + \varepsilon$ and 
$\mathcal S \notin \mathcal{SC}(L^2_{\oplus}(\RR_+))$; see \cite[Section 4.3]{benchimol1977feedback}.

\begin{proof}[Proof of Theorem~\ref{commutingCounter_intro}(i)]  
	
	We prove here that there exist two commuting $C_0$-semigroups
	$\mathcal T_1$ and $\mathcal T_2$ on a Hilbert space $H$ such that neither of them belongs to $\mathcal{SQC}(H)$ and $\mathcal T_1 \mathcal T_2 \in \mathcal{SC}(H)$.
	
	Let $\mathcal S = (S(t))_{t\geq0}$ be a quasi-contraction $C_0$-semigroup on a Hilbert space $\mathcal H$ satisfying $\sup_{t\geq0} \| S(t)\|<\infty$ and $\mathcal S \notin \mathcal{SC}(\mathcal H),$
	e.g. a Packel type semigroup or Benchimol's semigroup mentioned above. 
	Choose 
	$a \ge 0$
	such that $\|S(t)\| \leq e^{a t}$ for $t\geq0.$
		
	Now let $\mathcal T_1 = (T_1(t))_{t\geq0}$ and $\mathcal T_2 = (T_2(t))_{t\geq0}$ be the $C_0$-semigroups on $H := \oplus_{n=1}^\infty \left(\mathcal H \oplus \mathcal  H\right)$ given by 
	$$T_1(t) := \oplus_{n=1}^\infty \left(S(nt) \oplus e^{-ant}I\right),
	\qquad T_2(t):= \oplus_{n=1}^\infty \left(e^{-ant}I \oplus S(nt)\right), \qquad t\geq0.
	$$
	It is easy to see that $\mathcal T_1$ and $\mathcal T_2$ commute and that $\mathcal T_1 \mathcal T_2 \in \mathcal{SC}(H)$ since $e_{-a} \mathcal S$ is a contraction semigroup. As $\mathcal S \notin \mathcal{SC}(\mathcal H)$, using Chernoff's argument \cite[Section 3]{chernoff1976two} we infer that $\left(\oplus_{n=1}^\infty S(nt)\right)_{t\geq0} \notin \mathcal{SQC}\left(\oplus_{n=1}^\infty \mathcal H \right)$. Indeed, assume that there exist $b \in \mathbb R$ and an equivalent Hilbertian norm $\|\cdot\|_{\rm{eq}}$ on $\oplus_{n=1}^\infty \mathcal  H$ such that $\|\oplus_{n=1}^\infty e^{-b t}  S(nt)\|_{\rm{eq}} \leq 1$ for all $t>0$, and let $\|\cdot\|_{{\rm eq}, n}$ be the restriction of $\|\cdot\|_{\rm{eq}}$ to the $n$-th summand of  $\oplus_{n=1}^\infty \mathcal  H$. 
 Then the formula $\| h \|_{\mathscr H}^2 := \operatorname{LIM} \bigl[\| h \|_{{\rm eq}, n}^2 \bigr]$ for $h \in \mathcal H$, 
	defines an equivalent Hilbertian norm on $\mathcal H$ such that
	\begin{equation}\label{ChernoffEq}
		\|S(t) h\|_{\mathscr H}^2 = \operatorname{LIM} \bigl[\|e^{-\frac{b}{n} t} S(t) h\|_{{\rm eq}, n}^2\bigr] 
		\leq \operatorname{LIM} \bigl[\| h\|_{{\rm eq}, n}^2 \bigr] = \|h\|_{\mathscr H}^2, \qquad h \in \mathcal H, \, t \geq 0.
	\end{equation}
	This contradicts the assumption that $\mathcal S \notin \mathcal{SC}(\mathcal H)$.
	
	Since $\left(\oplus_{n=1}^\infty S(nt)\right)_{t\geq0}$ is the restriction of $\mathcal T_1$ to an invariant subspace, it then follows that $\mathcal T_1 \notin \mathcal{SQC}(H)$. (Otherwise, $\left(\oplus_{n=1}^\infty S(nt)\right)_{t\geq0}$ is a quasi-contraction semigroup in an appropriate norm.) Similarly, one concludes that $\mathcal T_2 \notin \mathcal{SQC}(H)$, which completes the proof.
\end{proof}

To show the failure of the ``if'' part of Theorem~\ref{split_finite}(ii) for commuting semigroups, we invoke the so-called generalized Bhat-Skeide interpolation introduced in 
\cite{bhat2015pure}
and developed, in particular, in 
\cite{dahya2024interpolation}.
To this aim, we interpolate Pisier's example \pcite[Theorem 1]{pisier1998joint} of two commuting operators 
on a Hilbert space, each similar to a contraction, such that their product is not similar to a contraction. 

Let us describe the generalized Bhat-Skeide interpolation in some more details. 
Given $m\in \NN$ and a family of commuting bounded operators $(T_k)_{k=1}^m$ on a Hilbert space $H$, 
the interpolation produces a family of commuting $C_0$-semigroups $(\mathcal T_k)_{k=1}^m$ 
on $L^2(\mathbb T^m)\otimes H$ satisfying
\begin{equation}\label{interpol_bs}
\prod_{k=1}^m T_k(n_k) = I  \otimes \left( \prod_{k=1}^m T_k^{n_k}\right), \qquad n_1, \ldots, n_m \in \ZZ_+.
\end{equation}
To construct $\mathcal T_k$, let $t\geq 0$ and $k \in \mathbb N_m$ be fixed.
 For every 
 \[
 \mathbf{z} = (z_1,\ldots, z_m) = (e^{2\pi i r_1},\ldots, e^{2\pi i r_m}) \in \mathbb T^m, \qquad r_1,\ldots, r_m \in [0,1),
\]
 define $\varphi_t^{(k)}: \mathbb T^m \to \mathbb T^m$ and  $p_t^{(k)}:\mathbb T^m \to \mathbb R$ by
\begin{align*}
	\left(\varphi_t^{(k)}(\mathbf{z})\right)_j &:= \begin{cases}
		e^{-2\pi i t} z_j, \quad &\mbox{if } j=k,
		\\ z_j, \quad &\mbox{if } j\neq k,
	\end{cases}
		\qquad \qquad  p_t^{(k)}(\mathbf{z}) = \chi_{[0, 1-\left\{t\right\} )} (r_k),
\end{align*}
where   $\left\{ t\right\}= t - \lfloor t \rfloor$.	Letting the operators $U_k(t), P_k(t) \in \linearOp(L^2(\mathbb T^m))$
be given by
\begin{align*}
	U_{k}(t) f = f \circ \varphi_t^{(k)}, \qquad P_{k}(t) f = p_t^{(k)} f, \qquad f \in L^2(\mathbb T^m),
\end{align*}
and observing that $U_{k}(n) = P_{k}(n) = I$ for 
for all $n \in \ZZ_+,$  we set
\begin{align}\label{BhatSkeideEq}
	T_k(t) &:= \left(U_{k}(t)  \otimes I\right) \left(P_{k}(t) \otimes T_k^{\lfloor t \rfloor}+ \left(I - P_{k}(t) \right)  \otimes T_k^{\lfloor t \rfloor +1}\right).
\end{align}
Then $\mathcal T_k$ is a  $C_0$-semigroup on $L^2(\mathbb T^m)\otimes  H,$
and the family $(\mathcal T_k)_{k=1}^m$  satisfies \eqref{interpol_bs}.
In this case, we will say that $(\mathcal T_k)_{k=1}^m$ interpolates $(T_k)_{k=1}^m.$
Moreover, for every $k \in \mathbb N_m$,  
the semigroup $\mathcal T_k$ is contractive  if and only if $T_k$ is a contraction.
For the proof of these properties, we refer the reader to \cite[Section 2]{dahya2024interpolation},
in particular to \cite[Lemma 1.4]{dahya2024interpolation} (see also \cite[Example 8.2]{shalit2023semigroups}).
Note that a similar argument shows that, in addition, $\mathcal T_k$ is bounded if and only if $T_k$ is power bounded,
i.e. $\sup_{n \ge 0}\|T_k^n\|<\infty$.

A nice feature of the Bhat-Skeide interpolation is that, for every $k \in \mathbb N_m,$ the similarity
of  $\mathcal T_k$ to a contraction semigroup is equivalent  to the similarity of $T_k$ to a contraction.
We formalize this fact in the next simple lemma,
which will be crucial for the proof of Theorem \ref{commutingCounter_intro}(ii).
\begin{lemma}\label{renorm_skeide}
Let $H$ be a Hilbert space, let $(T_k)_{k=1}^m \subset \mathcal L(H)$ be a family of commuting operators,
and let $(\mathcal T_k)_{k=1}^m\subset \mathcal L(L^2(\mathbb T)\otimes H)$ be the family of commuting semigroups given by
\eqref{BhatSkeideEq}.
\begin{itemize}
\item [(i)]  For every $k \in \mathbb N_m$ the semigroup $\mathcal T_k$ is similar to
a contraction semigroup on $L^2(\mathbb T^m) \otimes H$
if and only if $T_k$ is similar to a contraction on $H$. 
\item [(ii)] In addition, for any 
$(n_1, ..., n_m) \in \mathbb Z_+^m,$
$\prod_{k=1}^m T_k(n_k)$ is similar to a contraction if and only if $\prod_{k=1}^m T_k^{n_k}$ is so.
\end{itemize}
\end{lemma}
\begin{proof}
Let $k \in \mathbb N$ and $t \in \mathbb R_+$ be fixed, 
and let
 $R_k T_k R_k^{-1}$ be a contraction on $H$ for an invertible
$R_k.$ 
 Since $\left((I\otimes R_k) T_k(t) (I\otimes R_k^{-1})\right)_{t\geq0}$ interpolates $R_k T_k R_k^{-1},$
and the formula \eqref{BhatSkeideEq} preserves contractivity,
we infer that $\mathcal T_k$ is similar to a semigroup of contractions.
On the other hand, 
if $\mathcal T_k$ is similar to a contraction semigroup on $L^2(\mathbb T) \otimes H,$ then
 $\|T_k(t)\|_{{\rm eq}}\le 1$ for an equivalent norm $\|\cdot\|_{{\rm eq}}$ on $L^2(\mathbb T) \otimes H.$
Then it suffices to define the equivalent norm $\|\cdot\|_{\mathscr H}$ on $H$ as 
\begin{equation}\label{norm_h}
\|h\|_{\mathscr H}:=\|\chi_{\mathbb T}\otimes h\|_{{\rm eq}}, \qquad h \in H,
\end{equation}
and to note that 
\[
\|T_k h\|_{\mathscr H}=\|(I\otimes T_k)(\chi_{\mathbb T}\otimes h)\|_{{\rm eq}}
=\|T_k(1) (\chi_{\mathbb T}\otimes h)\|_{\rm{eq}}
\le \|\chi_{\mathbb T}\otimes h\|_{{\rm eq}}=\|h\|_{\mathscr H}
\]
for all $h \in H.$ 

For (ii), the claim follows directly from the interpolation property \eqref{interpol_bs} together with the elementary fact
 that for any Hilbert spaces $H_1, H_2$ and $T \in \mathcal L(H_2)$ 
the operator $I\otimes T$ on $H_1\otimes H_2$ is similar to a contraction  
if and only if $T$ possess this property (e.g. via an argument given above).
\end{proof}
We will also need the following version of \eqref{BhatSkeideEq} for $m=1$ not containing tensor products
and more convenient for some of our purposes.
Let $T \in \mathcal L(H)$ be fixed.
Identifying $L^2(\TT)\otimes H$ with $L^2(\TT,H)$ and omitting upper subscripts, 
we obtain the representations
\begin{align}\label{BSRep}
	(U(t)\otimes I) f = f\circ \varphi_t, \qquad
	(P(t) \otimes I) f = p_t f, \qquad
	((I\otimes T) f)(e^{2\pi i r}) = T (f (e^{2\pi i r})), 
\end{align}
for all $t\geq 0$,  $r\in [0,1)$ and for every simple function $f \in L^2(\TT,H)$. As the subspace of simple functions is dense in $L^2(\TT, H)$, one infers that \eqref{BSRep} holds for every $f\in L^2(\TT, H),$ and 
the corresponding interpolating semigroup $\mathcal T=(T(t))_{t \ge 0}$ is thus identified with
\begin{equation}\label{simple}
	\begin{aligned}
	(T(t) f) (e^{2\pi ir}) &= T^{\lfloor t \rfloor} ((p_t f)  (\varphi_t(e^{2\pi i r}))) + T^{\lfloor t \rfloor+1}(((I-p_t) f) (\varphi_t(e^{2\pi i r})))
	\\ &= \begin{cases}
		T^{\lfloor t \rfloor + 1} (f(e^{2\pi i(r-t)})), \qquad & r \in [0, \{t\}),
		\\T^{\lfloor t \rfloor} (f(e^{2\pi i(r-t)})), \qquad & r \in [\{t\}, 1),
	\end{cases}
	\end{aligned}
\end{equation}
where $f\in L^2(\TT, H)$ and $t\geq 0$.

Now we are ready to derive several implications of the Bhat-Skeide interpolation
in our context, and
we start with the proof of Theorem~\ref{commutingCounter_intro}(ii).

\begin{proof}[Proof of Theorem~\ref{commutingCounter_intro}(ii)] $ $ 
	
	We prove that there exist  two commuting $C_0$-semigroups $\mathcal S_1$ and $\mathcal S_2$ on
	on a Hilbert space $K$, both in $\mathcal{SC}(K)$, such that $\mathcal S_1 \mathcal S_2 \notin \mathcal{SQC}(K)$.
	
	Let $H :=\ell^2(\mathbb N^2).$
As shown by Pisier in \cite[Theorem 1]{pisier1998joint}, there exist commuting $T_1, T_2 \in \linearOp(H)$  
each  similar to a contraction, such that at the same time $T_1 T_2$ is not similar to a contraction.
	 Applying the generalized Bhat-Skeide interpolation to 
	$T_1$ and $T_2,$ we obtain two commuting $C_0$-semigroups on $L^2(\TT^2)\otimes H$, $\mathcal T_1 = (T_1(t))_{t\geq0}$ and $\mathcal T_2 = (T_2(t))_{t\geq0}$
	on $H$ satisfying \eqref{interpol_bs} for $m=2$.
	By Lemma \ref{renorm_skeide}(i) we have $\mathcal T_1, \mathcal T_2 \in \mathcal{SC}(L^2(\TT^2)\otimes H).$
	On the other hand, $T_1(1) T_2(1) = I  \otimes T_1 T_2,$ and by Lemma \ref{renorm_skeide}(ii),
		$T_1(1)T_2(1)$ is not similar to a contraction since $T_1 T_2$ is neither.
	Thus $\mathcal T_1 \mathcal T_2 = (T_1(t)T_2(t))_{t\geq0}$ is a $C_0$-semigroup on $L^2(\TT^2)\otimes H$ such that $\mathcal T_1 \mathcal T_2 \notin \mathcal{SC}\left(L^2(\TT^2)\otimes H\right).$ 
	Next, define another pair of commuting $C_0$-semigroups $\mathcal S_1 = (S_1(t))_{t\geq0}$ and $\mathcal S_2 =(S_2(t))_{t\geq0}$ on $K := \oplus_{n=1}^\infty \left(L^2(\mathbb T^2)  \otimes H \right)$ as
	$$
	S_1(t) := \oplus_{n=1}^\infty T_1(nt), \qquad S_2(t) := \oplus_{n=1}^\infty T_2(nt), \qquad t \geq 0.
	$$
	In view of $\mathcal T_1, \mathcal T_2 \in \mathcal{SC}(L^2(\TT^2)\otimes H),$ it follows that $\mathcal S_1, \mathcal S_2 \in \mathcal{SC}(K)$ . Moreover, since $\mathcal T_1 \mathcal T_2$ does not belong to $\mathcal{SC}\left(L^2 (\TT^2) \otimes H\right)$, using Chernoff's direct sum construction as in the proof of Theorem~\ref{commutingCounter_intro}(i), see \eqref{ChernoffEq}, we conclude that $\mathcal S_1 \mathcal S_2 \notin\mathcal{SQC}(K)$, thus completing the proof.
\end{proof}

\begin{remark}
	It is worth mentioning that, using Pisier's example, \cite[Example 4.3]{badea2005operators} provides a pair of bounded operators $T_1$ and $T_2$ on a  Hilbert space satisfying 
	$T_1 T_2 T_1=T_1$ and $T_2=T_2 T_1 T_2$, both of them similar to contractions, yet not simultaneously similar to contractions. We omit the discussion of a semigroup counterpart of this example.
\end{remark}

\subsection{Bhat-Skeide interpolation, quasisimilarity and a Peller's question}\label{MTSubsect}

Tensor products provide a natural way to transfer intricate constructions from the setting of single operators to that of 
$C_0$-semigroups. In this section, we adapt several well-known examples of operators which, despite having additional properties, are not similar to contractions. Using the Bhat–Skeide interpolation introduced above, these single operator examples yield new examples in the semigroup framework. This method contrasts with the traditional approach, which relies on ad hoc adaptations of discrete techniques, often technically demanding, so that many continuous counterparts have been absent from the literature.

We begin with the continuous analogue of a result in \cite{muller2007quasisimilarity}, showing
the existence of power bounded operators that are not \textit{quasi-similar} to contractions.  
Recall that given Hilbert spaces $H$ and $K$ and operators $T \in \linearOp(H)$ and $S \in \linearOp(K)$, we write $S \prec T$ if there exists a bounded, injective operator $R \in \linearOp(K,H)$ with dense range such that $TR = RS$. When both $S \prec T$ and $T \prec S$ hold, $T$ and $S$ are said to be quasi-similar. Clearly, if two operators are similar to each other, they are also quasi-similar, but the converse in general does not hold.
To show that the two operators are not quasi-similar it is enough to find a property stable under quasi-similarity, and possessed by only one operator from the pair. As shown in \cite[Theorem 1.1 \& Lemma 3.3]{muller2007quasisimilarity}, the so-called Blum-Hanson property can serve for this purpose.  For $T \in \mathcal L(H)$, one of the formulations of the Blum-Hanson property requires that for every $h \in H$
the weak convergence of $T^n h$ implies the strong convergence of $\frac{1}{m} \sum_{k=1}^m T^{n_k}h$
for all (increasing) subsequences $(n_k)_{k=1}^\infty \subset \mathbb N.$
In \cite[Example 2.1]{muller2007quasisimilarity} it was constructed an operator $M$ on a 
Hilbert space $H$
such that $\lim_{n\to \infty} T^n = 0$ in the weak operator topology, $\inf_{n\in \NN} \|T^n h\| > 0$ for each $h\in H\setminus\{0\}$, and 
$\inf_{m \in \mathbb N} \|\frac{1}{m} \sum_{k=1}^m T^{\widehat n_k}h\|>0$ for some $(\widehat n_k)_{k=1}^\infty \subset \NN$ and $h \in H,$
so that $T$ does not satisfy the Blum-Hanson property.
 Since Hilbert space contractions possess the Blum-Hanson property, $T$  is then not quasi-similar to a contraction on $H,$  and it will be a building block for the subsequent argument.

\begin{proof}[Proof of Theorem~\ref{MTTh1}]
	We show that there exist a $C_0$-semigroup $\mathcal T = (T(t))_{t\geq0}$ on a Hilbert space $H$ such that $\lim_{t\to \infty} T(t)=0$ in the weak operator topology, $\inf_{t>0} \|T(t)h\| > 0$ for each $h\in H\setminus\{0\}$, and $S(t) \prec T(t)$, $t>0$, holds for no contraction $C_0$-semigroup $\mathcal S = (S(t))_{t\geq0}$.
	
Let $\mathcal H$ and $T \in \linearOp(\mathcal H)$ denote, respectively, the Hilbert space and the operator from \cite[Example 2.1]{muller2007quasisimilarity} described above.
Applying the Bhat-Skeide interpolation to $T$, construct the corresponding Bhat-Skeide interpolating $C_0$-semigroup $\mathcal T = (T(t))_{t\geq0} \subset \linearOp (H)$, where $H := L^2(\TT) 
\otimes \mathcal H.$ Since $(T^n)_{n \ge 0}$
converges to zero in the weak operator topology it follows from 
\eqref{simple} 
	that $(T(t)f)_{t \ge 0}$ converges to zero in the weak operator topology for all simple functions
$f \in H,$
and then on the whole of $H$ by density arguments.
	Moreover, by \eqref{simple} it is easy to see that
	\begin{align*}
		\inf_{n\in \NN} \|T(n)f\| = \inf_{n\in \NN} \|(I \otimes T^n) f\|>0, \qquad f \in H \setminus\{0\}. 
	\end{align*}
	In fact, assuming otherwise, we would obtain that $\lim_{k\to \infty} T^{n_k} f(z)=0$ for a.e. $z\in \TT$ and for a subsequence $(n_k)_{k \in \NN} \subseteq \NN$, which contradicts the assumption that $\inf_{n\in \NN} \|T^n h\|>0$ for all $h\in H \setminus\{0\}$. As $\sup_{t> 0} \|T(t)\| < \infty,$ 
	we have that
	$$\inf_{t>0} \|T(t)f\|>0, \qquad f \in H \setminus\{0\}.
	$$	
	Similarly, in view of \eqref{simple}, the sequence $\frac{1}{m} \sum_{k=1}^m T(\widehat n_k) = I \otimes \left( \frac{1}{m} \sum_{k=1}^m  T^{\widehat n_k} \right)$ does not converge in the strong operator topology as $m\to \infty$ due to the choice of $T$. Consequently, $T(1)$ does not satisfy the Blum-Hanson property, and it follows from \cite[Theorem 1.1 \& Lemma 3.3]{muller2007quasisimilarity} that there is no contraction  $C$ on $H$ such that $C \prec T(1)$, completing the proof.
\end{proof}
An example of $\mathcal T$ which cannot be intertwinned with a contraction even in two possible ways
can be constructed using \cite[Example 3.7]{muller2007quasisimilarity}. We omit the argument, which is similar
to the proof of Theorem \ref{MTTh1}.

\begin{remark}
	Let $\varepsilon > 0$ be given. With a minor modification to the proof of Theorem~\ref{MTTh1}, one may construct a $C_0$-semigroup $\mathcal T = (T(t))_{t\geq0}$ 
	which apart from the properties stated in the theorem satisfies 
	$$
	\sup_{t\geq0} \|T(t)\| \leq 1 + \varepsilon.
	$$
		This can be achieved by 
	applying the Bhat-Skeide interpolation to  $T,$ which, in addition to the properties in 
	\cite[Example 2.1]{muller2007quasisimilarity}, satisfies $\sup_{n \in \NN} \|T^n\| \leq 1 + \varepsilon$. In fact, developing a simplified version of \cite[Example 2.1]{muller2007quasisimilarity}, \cite[Example 5.12]{van2009blum} offers an operator $T$ on a Hilbert space $\mathcal H$ such that $\sup_{n \in \NN} \|T^n\| \leq 1 + \varepsilon$, $\lim_{n\to \infty} T^n = 0$ in the weak operator topology, and that does not have the Blum-Hanson property. Applying the arguments from \cite{van2009blum} to the original example in \cite[Example 2.1]{muller2007quasisimilarity}, one produces our operator $T$. Alternatively, one can easily modify \cite[Example 2.1]{muller2007quasisimilarity}
	for this purpose.

\end{remark}

Let \(T\) be a power bounded operator on a Hilbert space \(H\) such that
\begin{equation}\label{sepzeroEq}
	\inf_{n\in \NN} \|T^n h\| > 0 \quad \text{and} \quad \inf_{n\in \NN} \|(T^\ast)^n h\| > 0, \qquad h \in H\setminus\{0\}.
\end{equation}
By a classical result due to Sz.-Nagy and Foia\c{s}, $T$ is quasi-similar to a unitary operator. At the same time, providing a negative answer to a question by K\'erchy \cite[Question 1]{kerchy1989isometric}, it was shown in \cite[Example 4.1]{muller2007quasisimilarity} that there exists a power bounded operator $T$ on $H$ not similar to a contraction and satisfying \eqref{sepzeroEq}.  We now establish the semigroup analogue of such a result.

\begin{proposition}\label{MTProp2}
	There exists a $C_0$-semigroup $\mathcal T = (T(t))_{t\geq0}$ on  a Hilbert space $H$ satisfying $\sup_{t\ge 0} \|T(t)\| < \infty$, $\inf_{t>0} \|T(t)h\|>0$ and $\inf_{t>0} \|T(t)^\ast h\| > 0$ for each $h\in H\setminus\{0\}$, and such that $\mathcal T \notin \mathcal{SC}(H)$. 
\end{proposition}
\begin{proof}
	Let $T$ be the power bounded operator on a Hilbert space $\mathcal H$ given in 
	\cite[Example 4.1]{muller2007quasisimilarity}, so that $T$ is not similar to a contraction
	and satisfies \eqref{sepzeroEq}.
	Applying as above the Bhat-Skeide interpolation to $T$ we construct the $C_0$-semigroup $\mathcal T  =(T(t))_{t\geq0} \subset \linearOp(H)$, with $H := L^2(\TT) \otimes \mathcal H$. 
	Arguing as in the proof of Theorem~\ref{MTTh1}, we infer that $\sup_{t>0} \|T(t)\| < \infty$, $\inf_{t>0} \|T(t)h\|>0$ and $\inf_{t>0} \|T(t)^\ast h\| > 0$ for each $h\in H\setminus\{0\}$. Since $T(1) = I \otimes T$, from Lemma \ref{renorm_skeide}(i) it follows that $\mathcal T \notin \mathcal{SC}(H)$, completing the proof. 
\end{proof}

Developing further the interpolation approach, recall that Peller asked in \cite{peller1982estimates} whether, for every \(\varepsilon > 0\) and any power bounded operator \(T\) on a Hilbert space $H$, there exists an equivalent Hilbertian norm \(\|\cdot\|_{\rm{eq}}\) on \(H\) such that  
\[
\sup_{n\in \mathbb{N}} \|T^n\|_{\rm{eq}} \leq 1 + \varepsilon.
\]  
This long-standing problem was solved by Kalton and Le Merdy  in \cite{kalton2002solution} by constructing a counterexample. Using their counterexample and the Bhat-Skeide interpolation,
we show that the Kalton-Le Merdy construction transfers to the semigroup framework
and a similar phenomena holds for $C_0$-semigroups.

\begin{proposition}\label{KaltonLeMProp}
	Let $\varepsilon>0$ be fixed. Then there exist a $C_0$-semigroup $\mathcal T = (T(t))_{t\geq0}$
	on a Hilbert space $H$ such that $\sup_{t>0} \|T(t)\| < \infty$ and  
	$$\sup_{t>0} \|T(t)\|_{\rm{eq}} > 1 +\varepsilon,
	$$
	for every equivalent Hilbertian norm $\|\cdot\|_{\rm{eq}}$ on $H$.
\end{proposition}
\begin{proof}
	Given $\varepsilon>0,$ we infer from \cite[Theorem 3.4]{kalton2002solution} that there are a Hilbert space $\mathcal H$ and $T\in \mathcal{L}(\mathcal H)$  such that $T$ is power bounded 
	and
	$\sup_{n \in \NN} |||T^n||| > 1 +\varepsilon$ for every equivalent Hilbertian norm $|||\cdot |||$ on $\mathcal H$. Construct the $C_0$-semigroup $\mathcal T =(T(t))_{t\geq0} \subset \linearOp(H)$, where $H:= L^2(\TT) \otimes \mathcal H$, by applying the Bhat-Skeide interpolation to $T$.
	Since $T$ is power bounded, we have $\sup_{t\geq0} \|T(t)\| < \infty.$
	Moreover, given an equivalent norm $\|\cdot\|_{{\rm eq}}$ on $H$ we argue as in the proof of Lemma \ref{renorm_skeide} and define the  equivalent Hilbertian norm $|||\cdot|||$ on $\mathcal H$ 
	by $|||h||| := \| \chi_\TT \otimes h\|_{\rm{eq}}$, $h \in \mathcal H$.
	Then   
	$$\sup_{n \in \NN} \|T(n)\|_{\rm{eq}} = \sup_{n\in \NN} \| I \otimes T^n\|_{\rm{eq}} \geq \sup_{n\in \NN} 
	|||T^n||| > 1+\varepsilon,
	$$
	as required.
\end{proof}

\subsection{(Quasi-)nilpotent semigroups not similar to contraction ones}\label{counterexamplesSub}

The problem of similarity of $C_0$-semigroups to semigroups of contractions differs substantially from the corresponding problem for individual operators.
To support this claim, we present two examples of semigroups whose properties might suggest similarity to contraction semigroups, but which in fact do not even ensure similarity to quasi-contraction semigroups.
Specifically, using Theorem~\ref{split_finite}, we prove that neither quasi-nilpotency nor nilpotency is sufficient for a semigroup to be similar to a contraction one, even under additional assumptions such as immediate compactness and holomorphicity in $\CC_+$ (the latter in the quasi-nilpotent case).
This is in sharp contrast with the discrete setting, where compactness and/or (quasi-)nilpotency 
ensure similarity to a contraction by a direct argument. Moreover, since only a (quasi~\!\!-~\!\!)\!\!~nilpotent operator $T$ leads to a (quasi\!\!~-\!\!~)\!\!~nilpotent semigroup $\mathcal T$ through \eqref{simple}, the examples we present here cannot be obtained via the Bhat-Skeide interpolation.

Recall first some basic semigroup terminology crucial for the sequel. A $C_0$-semigroup $\mathcal T = (T(t))_{t\geq0}$ is said to be 
\begin{enumerate}
	\item [(i)] \textit{quasi-nilpotent} if $\omega_0(\mathcal T) = -\infty$, or equivalently, $\sigma(T(t)) = \{0\}$ for some (hence all) $t>0$;
	\item [(ii)] \textit{nilpotent} if there exists $\tau>0$ such that $T(t) = 0$ for all $t\geq \tau$;
	\item [(iii)] \textit{immediately compact} if $T(t)$ is a compact operator for all $t >0$;
	\item [(iv)]  \textit{bounded holomorphic of angle} $\theta \in (0,\pi/2]$ if $\mathcal T$ has a holomorphic extension 
	to the sector 
	$\Sigma_{\theta} := \{z \in \CC \setminus\{0\} \, : \, |\operatorname{arg} z| < \theta\}$
	such that
	\begin{align*}
		\sup \{\|T(z)\| \, : \, z \in \overline{\Sigma_{\theta'}}\} < \infty \qquad \mbox{for all } \theta' \in (0,\theta).
	\end{align*}
\end{enumerate}

It is easy to prove that the above properties are preserved under tensor products in the following sense. 
\begin{lemma}\label{NilpotCompacLem}
	Let $\mathcal T = (T(t))_{t\geq0}$ and $\mathcal S = (S(t))_{t\geq0}$ be $C_0$-semigroups on Hilbert spaces $H$ and $K$ respectively. Then the following holds.
	\begin{enumerate}
		\item [(i)] $\mathcal T  \otimes \mathcal S$ is (quasi-)nilpotent if and only if either $\mathcal T$ is (quasi\nobreakdash-)nilpotent or $\mathcal S$ is (quasi\nobreakdash-)nilpotent.
		\item [(ii)] $\mathcal T  \otimes \mathcal S$ is immediately compact if and only if both $\mathcal T$ and $\mathcal S$ are immediately compact.
		\item [(iii)] Fix $\theta \in (0,\pi/2]$. Then $\mathcal T  \otimes \mathcal S$ is bounded holomorphic of angle $\theta$ if both $\mathcal T$ and $\mathcal S$ are bounded holomorphic of angle $\theta$.
	\end{enumerate}
\end{lemma}
\begin{proof}	
	The claim in (i) concerning nilpotency is trivial, and the statement addressing quasi-nilpotency is a straightforward consequence of the identity $\sigma(T(1) \otimes S(1)) = \sigma(T(1)) \sigma(S(1))$, valid by \eqref{prod_spect}.

	The property (ii) follows from Holub \cite[Proposition 3.1]{holub1970tensor} or Apiola \cite[Theorem 4.5]{apiola1973tensorproduct} in the more general setting of Banach spaces. We give a short and illustrative proof of this property for operators on Hilbert spaces. For the ``if'' part, if $(A_n)_{n=1}^\infty \subset \linearOp(H)$ and $(B_n)_{n=1}^\infty \subset \linearOp(K)$ are sequences of finite-rank operators such that $\lim_{n\to \infty} A_n = A$ and $\lim_{n\to \infty} B_n =B$ in the uniform operator topology, then 
	$(A_n \otimes B_n)_{n=1}^\infty \in \linearOp(H\otimes K)$ is a sequence of finite-rank operators\footnote{In fact, let $T_{x,y}$ be the rank-one operator given by $z\mapsto \langle z,x\rangle y$, where $x,y$ are two fixed vectors. Then $T_{x_1,y_1}  \otimes T_{x_2,y_2}= T_{x_1 \otimes x_2, y_1\otimes y_2}$ is also rank-one.} such that $\lim_{n\to\infty}  A_n\otimes B_n = A\otimes B$ in the uniform operator topology. For the ``only if'' part, 
	fix $ A \in \linearOp(H)\setminus\{0\}$, $B \in \linearOp(K)\setminus\{0\}$ and $k_0 \in K$ satisfying $Bk_0 \neq 0.$ Define the operators $Q\in \linearOp(H, H \otimes K)$ and  $P \in \linearOp( H  \otimes K, H)$ by
	setting
	$$Qh = h \otimes k_0, \quad h \in H,   
	\qquad P\left(h \otimes k \right) =  \langle k, Bk_0\rangle h,  \quad h \in H, \,\, k \in K,
	$$
 extending $P$ to $H \otimes K$ by linearity and density, and denoting the extension by the same symbol.	It is readily seen that $P (A\otimes B) Q$ is a (nonzero) scalar multiple of $A$, thus 
$A$ is a compact operator.  The same argument shows that $B$ is also compact. 
	
	Finally, to show (iii),
		it suffices to note that $\mathcal T \otimes \mathcal I$ and $\mathcal I  \otimes \mathcal S$ are two commuting bounded holomorphic $C_0$-semigroups of angle $\theta$ on $H  \otimes K$ (their holomorphic extensions are given by $(T(z)\otimes I)_{z \in \Sigma_\theta}$ and $(I  \otimes S(z))_{z\in \Sigma_\theta}$ respectively).
Since $\mathcal T \otimes \mathcal S =(\mathcal T \otimes \mathcal I) (\mathcal I \otimes \mathcal S),$
the statement follows.
\end{proof}

For further considerations, we will need a pair of auxiliary semigroups.

\begin{example}\label{examplesSemi} $ $ \newline

	\begin{enumerate}
		\item Let $\mathcal T_{RL} = (T_{RL}(t))_{t\geq0}$ be the Riemann-Liouville $C_0$-semigroup on $L^2[0,1]$, given by
		$$(T_{RL}(t)f)(x) = \frac{1}{\Gamma(t)} \int_0^x (x-y)^{t-1} f(y)\, dy, \qquad x \in [0,1], \, f \in L^2[0,1], \, t >0. 
		$$
		It is well known that $\mathcal T_{RL}$ is immediately compact, quasi-nilpotent and bounded holomorphic of angle $\frac{\pi}{2},$ 
				see for instance \cite[Theorems 3.1 and 3.4]{carvalho2022riemann} (containing, in fact, more general results).
		\item  Let $\mathcal H$ be a separable Hilbert space and let $(h_n)_{n=1}^\infty$ be a conditional basis of $\mathcal H$. For every $t \ge 0$ define a bounded operator $T_{LM}(t)$ on $\mathcal H$ by setting 
		\begin{equation}\label{leMerdySem}
			T_{LM}(t) h_n: = e^{-2^n t} h_n, \qquad n\in \NN, \, t \geq 0,
		\end{equation}
		extending this map to $\mathcal H$ by linearity and density, and denoting the extension by the same symbol. Then, as Le Merdy proved in \cite{lemerdy2000bounded}, $\mathcal T_{LM} = (T_{LM}(t))_{t\geq0}$ is an  immediately compact, exponentially stable $C_0$-semigroup on $\mathcal H$ such that $\mathcal T_{LM} \notin \mathcal{SC}(\mathcal H)$.
			In addition, it is easy to see that such a semigroup is sectorially bounded holomorphic of angle $\frac{\pi}{2}$; see, for instance \cite[Proposition 3.5]{fackler2015regularity}. 
		Note that moreover, by Proposition~\ref{contractiveStableProp}, it follows that $\mathcal T_{LM} \notin \mathcal{SQC}(\mathcal H)$. Alternatively, to show that $\mathcal T_{LM} \notin \mathcal{SQC}(\mathcal H)$,  one may apply the arguments given in \cite{lemerdy2000bounded} for $\mathcal T_{LM}$ to prove that $e_{d} \mathcal T_{LM}  \notin \mathcal{SC}(\mathcal H)$ for each $d \in \RR.$ 
	\end{enumerate}
\end{example}

Given Theorem \ref{split_finite} and Lemma \ref{NilpotCompacLem}, and having in mind semigroups $\mathcal T_{RL},$ $\mathcal T_{LM},$ we are now ready to prove Theorem \ref{Counter_intro}, whose proof becomes comparatively direct.

\begin{proof}[Proof of Theorem~\ref{Counter_intro}] $ $ 
	
	We show first that there exist a Hilbert space $H$ and a nilpotent, immediately compact $C_0$-semigroup $\mathcal T$ on $H$ such that $\mathcal T \notin \mathcal{SQC}(H)$.
	
	First, note that the Riemann-Liouville semigroup $\mathcal T_{RL}$ commutes with the right shift semigroup $\mathcal R$ on $L^2[0,1]$, so $\mathcal T_{RL} \mathcal R = (T_{RL}(t)R(t))_{t\geq0}$ 
	is a $C_0$-semigroup on $L^2[0,1]$. Moreover, as $\mathcal R$ is nilpotent and $\mathcal T_{RL}$ is immediately compact,  $\mathcal T_{RL} \mathcal R$ is nilpotent and immediately compact. Consider the $C_0$-semigroup $\mathcal T := (\mathcal T_{RL} \mathcal R)  \otimes \mathcal T_{LM}$ on $H := L^2[0,1]\otimes \mathcal H$, where  $\mathcal T_{LM}$ is given by \eqref{leMerdySem}. Since $\mathcal T_{RL}\mathcal R$ is nilpotent, from Lemma~\ref{NilpotCompacLem}(i)  it follows that $\mathcal T$ is nilpotent as well, and since both $\mathcal T_{RL} \mathcal R$ and $\mathcal T_{LM}$ are immediately compact,
Lemma~\ref{NilpotCompacLem}(ii)	implies that $\mathcal T$ is so.
On the other hand, as $\mathcal T_{LM} \notin \mathcal{SQC}(\mathcal H)$, Theorem~\ref{split_finite}(i) then yields that $\mathcal T \notin \mathcal{SQC}(H)$. \medskip
		
	We prove now that there exists a quasi-nilpotent, immediately compact and bounded holomorphic of angle $\pi/2$ $C_0$-semigroup $\mathcal S$ on a Hilbert space $H$ such that $\mathcal S \notin \mathcal{SQC}(H)$.
		Define $\mathcal S := \mathcal T_{RL}  \otimes \mathcal T_{LM}$, where $\mathcal T_{RL}$ is the Riemann-Liouville semigroup on $L^2[0,1]$ and $ \mathcal T_{LM}$ is the semigroup given by \eqref{leMerdySem}. One has by Lemma~\ref{NilpotCompacLem} that $\mathcal S$ is quasi-nilpotent since $\mathcal T_{RL}$ is quasi-nilpotent, and that
 $\mathcal S$ is immediately compact and bounded holomorphic of angle $\pi/2$ since both $\mathcal T_{RL}$ and $\mathcal T_{LM}$ are so. On the other hand, since $\mathcal T_{LM} \notin \mathcal{SQC}(\mathcal H)$, Theorem~\ref{split_finite}(i) shows that $\mathcal S \notin \mathcal{SQC}(H)$, completing the proof.
\end{proof}

Our approach leads to new (counter-)examples also in the discrete setting. We restrict ourselves by a sample.
The following question was posed in \cite[Problem 3]{vu1997almost}: Let $A$ be a power bounded operator on a Hilbert space $H$ such that $A$ commutes with an injective and compact operator $B$. Is $A$ similar to a contraction?

We show below that the cogenerator of the $C_0$-semigroup $\mathcal{S}$ constructed in Theorem~\ref{Counter_intro} provides a negative answer, even under the additional assumptions that $B$ is quasi-nilpotent and belongs to the bicommutant $\{A\}''$ of $A$. 
Note that if the quasi-nilpotency assumption is dropped, then the cogenerator of the semigroup defined by \eqref{leMerdySem} would suffice for that purpose.

\begin{proposition}
	There exist bounded linear operators $A$ and $B$  on a Hilbert space $H$ satisfying the following.
	\begin{itemize}
		\item [(i)] $B$ is compact, injective and quasi-nilpotent.
		\item [(ii)] $B \in \{A\}''$.
		\item [(iii)] $A$ is power bounded and not similar to a contraction.
	\end{itemize}
\end{proposition}  
\begin{proof}
	Let $\mathcal S = (S(t))_{t\geq0}$ be the $C_0$-semigroup constructed in Theorem~\ref{Counter_intro}, and let $\mathcal A$ be its generator. As $\mathcal S$ is a bounded semigroup that is not in $\mathcal{SC}(H)$, its cogenerator $A$, defined by
	\begin{equation}\label{cogenEq}
		A = (\mathcal A+I) (\mathcal A-I)^{-1} = I + 2 (\mathcal A - I)^{-1}, 
	\end{equation}
	is not similar to a contraction; see, for instance, \cite[p. 109]{gomilko2017inverses}. 
	As $\mathcal S$ is holomorphic and sectorially bounded, one infers that $A$ is power bounded, see e.g. \cite[Theorem 5]{gomilko2004cayley}. 
	
	Fix $\tau>0$ and set $B = S(\tau)$, so that $B$ is compact and quasi-nilpotent. Let $C$ be a bounded operator on $H$ that commutes with $A$. It follows from \eqref{cogenEq} that $C$ commutes with $(I-\mathcal A )^{-1},$  and then with $(\lambda -\mathcal A)^{-1}$ for every $\lambda$ in the resolvent set of $\mathcal A$; see, for example, \cite[Proposition B.7]{arendt2011vector}. Since by the Post-Widder inversion formula for $\mathcal S$ (see, for instance, \cite[Corollary 3.3.6]{arendt2011vector}), 
	$$
	S(t)x=\lim_{n \to \infty}(I-t\mathcal A/n)^{-n}x, \qquad x \in H, \quad t \ge 0,
	$$
we conclude that $C$ commutes with $\mathcal S$, and thus  $B \in  \{A\}''$. 
	
	Finally, note that $B = S(\tau) = T_{RL}(\tau)  \otimes T_{LM}(\tau) = (T_{RL} (\tau)  \otimes I) (I  \otimes T_{LM}(\tau))$, where $\mathcal T_{RL} = (T_{RL}(t))_{\geq0}$ is the Riemann-Liouville semigroup, 
	and $\mathcal T_{LM} = (T_{LM}(t))_{t\geq0}$ with $T_{LM}(t)$ given by \eqref{leMerdySem}.  Since both $T_{RL}(\tau)$ and $T_{LM}(\tau)$ are injective, $B$ is injective as well,  and the claim follows.
\end{proof}

\section{Similarity and finite tensor products of operators}

In this section we adapt the methods of the preceding sections to obtain counterparts of our main results for tensor products of operators. This case is simpler than that of semigroups and, in fact, part of it already appears implicitly in \cite[Theorem 2.3]{paulsen1995centered}. Our approach, however, differs from that in \cite{paulsen1995centered} and yields explicit equivalent norms arising from similarities, which makes it natural to develop a unified treatment for both settings.

We give below the counterpart of one of our main results (Theorem~\ref{split_finite}(ii)) for the tensor product of a finite number of bounded operators.
The argument is similar to the proof of Theorem \ref{split_finite}, but substantially simpler,
and it does not require auxiliary results such as Proposition \ref{contractiveStableProp}. To put it in contrast to \cite{paulsen1995centered} we provide some details.

We now state the analogue of Theorem~\ref{split_finite}(ii) for finite tensor products of bounded operators. The proof follows the same general strategy as that of Theorem~\ref{split_finite}, but is considerably simpler and does not rely on auxiliary tools such as Proposition~\ref{contractiveStableProp}. To emphasize the difference with \cite{paulsen1995centered}, we include some details.
\begin{theorem}\label{Discrete1}
	For each $k \in \mathbb N_m$, let $H_k$ be a Hilbert space and let $T_k \in \linearOp(H_k)$. 
	Then $ \otimes_{k=1}^m T_k$ is similar to a contraction on $ \otimes_{k=1}^{m}H_k$ if and only if there exists 
	$(\alpha_k)_{k=1}^m \subset (0,\infty)$ with $\prod_{k=1}^m \alpha_k =1$ such that $\alpha_k T_k$ is similar to a contraction on $H_k$ for every $k \in \mathbb N_m$. 
\end{theorem}

\begin{proof}
To prove the ``if'' implication of the theorem,	let $(\alpha_k)_{k=1}^m \subset (0,\infty)$ be such that $\prod_{k=1}^m \alpha_k = 1$  and $\alpha_k T_k$ is similar to a contraction on $H_k$ for each $k\in \mathbb N_m$. Then, in view of Lemma~\ref{moreEquivalentTensorRemark}, we can assume that 
$\|\alpha_k T_k\| \leq 1$ for all $k \in \mathbb N_m.$ This yields
	$$\left\| \otimes_{k=1}^m T_k \right\| = \left\| \otimes_{k=1}^m \alpha_k T_k \right\| = \prod_{k=1}^m \|\alpha_k T_k\| \leq 1,
	$$
as required.
	
	We now turn to the “only if” implication, starting with $m=2$. Let $H_1$ and $H_2$ be two Hilbert spaces, and let $T_1 \in \linearOp(H_1)$ and $T_2 \in \linearOp(H_2)$ be such that $T_1  \otimes T_2$ is similar to a contraction. 
By \eqref{prod_spect},  we have	$r(T_1) r(T_2) = r(T_1  \otimes T_2) \leq 1$. If $r(T_1) r(T_2) < 1$, then  replacing $T_1$ and $T_2$ by $\alpha T_1$ and  ${\alpha}^{-1}T_2$ with  $\alpha \in (r(T_2), r(T_1)^{-1}),$ respectively, we note  that $r(\alpha T_1) < 1$ and $ r(\alpha^{-1} T_2) < 1$, so both $\alpha T_1$ and $\alpha^{-1} T_2$ are similar to contractions by Rota's theorem.
	If $r(T_1) r(T_2) = 1$, we may assume without loss of generality that $r(T_1)= r(T_2) = 1$ by replacing $T_1$ with $\frac{T_1}{r(T_1)}$ and $T_2$ with $\frac{T_2}{r(T_2)}.$ Hence, $\|T_1^n\| \geq1$ and $\|T_2^n\|\geq 1$ for all $n\in \NN$, and since  $T_1 \otimes T_2$ is power bounded by assumption,
	$$\sup_{n\in \NN} \|T_1^n\| = \sup_{n\in \NN} \left(\frac{\| (T_1 \otimes T_2)^n \|}{\| T_2^n \|} \right) \leq \sup_{n\in \NN} \| (T_1 \otimes T_2)^n \| < \infty.
	$$

	Let $\|\cdot\|_{\rm{eq}}$ be a Hilbertian norm on $H_1  \otimes H_2$ such that $\|T_1  \otimes T_2\|_{\rm{eq}} \leq 1$ and
	$$ \|\cdot\| \leq \|\cdot\|_{\rm{eq}} \leq \mathcal C \|\cdot\|,
	$$
	where $\mathcal C = \mathcal C(T_1\otimes T_2)$ is the similarity constant of $T_1 \otimes T_2$.
	Since $r(T_1)=1$, by Lemma \ref{vanNerveenLemma}(ii), for every $n \in \mathbb N$ there exists $h_n \in  H$ such that $\|h_n\| = 1$ and $\|T_1^k h_n\| \geq 1/2$ for all $k \in \mathbb N_n$. Now, for each $n \in \NN$, define the inner product $\langle \cdot, \cdot\rangle_n$
	on $H_2$ by
	\begin{align*}
		\langle f, g\rangle_n := \frac{1}{n+1} \sum_{k=0}^n \langle  T_1^k h_n \otimes f,  T_1^k h_n \otimes g \rangle_{\rm{eq}}, \qquad f, g \in H_2.
	\end{align*} 
	Then by the choice of $(h_n)_{n=1}^\infty$
	we have
		$$D := \sup_{n\in\NN} \frac{1}{n+1} \sum_{k=0}^n \|T_1^k h_n\|^2< \infty \qquad \text{and} \qquad  d := \inf_{n\in\NN} \frac{1}{n+1} \sum_{k=0}^n \|T_1^k h_n\|^2 >0.
	$$
	Hence,
	\begin{align}\label{SingleOpEq1}
		\|f\|_n^2 &= \frac{1}{n+1} \sum_{k=0}^n \|T_1^k h_n \otimes f\|_{\rm{eq}}^2 
		\leq  \mathcal C^2 D \|f\|^2, \qquad f \in H_2, \, n \in \NN,
	\end{align}
	and similarly,
	\begin{equation}\label{SingleOpEq2}
		\|f\|_n^2 \geq d \|f\|^2, \qquad f \in H_2, \, \, n \in \NN.
	\end{equation}
	 Moreover, arguing as in the proof of Theorem \ref{split_finite}(ii),
	\begin{align*}
		\|T_2f\|_n^2 
		 &= \frac{1}{n+1} \|h_n \otimes T_2 f\|_{\rm{eq}}^2 + \frac{1}{n+1}\sum_{k=1}^n \|T_1^k h_n \otimes T_2f\|_{\rm{eq}}^2\\
		\\ &\leq \left(1 + \frac{\mathcal C^2 \|T_2\|}{ (n+1) d}\right) \|f\|_n^2, \qquad f \in H_2, \, n \in \NN.
	\end{align*}
	So, setting   
	$$\|f\|_{\mathscr{H}_2}^2 := \operatorname{LIM} \left[ \|f\|_n^2\right], \qquad f \in H_2,
	$$
	and using \eqref{SingleOpEq1} along with \eqref{SingleOpEq2}, we infer that $\|\cdot\|_{\mathscr H_2}$ is a Hilbertian norm on $H_2$ satisfying
	\begin{align*}
		\|T_2 f\|_{\mathscr{H}_2}^2 
		\leq \operatorname{LIM} \left[ \left(1 + \frac{\mathcal C^2 \|T_2\|}{ (n+1) d} \right) \|f\|_n^2\right]
		= \|f\|^2_{\mathscr{H}_2}, \qquad f \in H_2.
	\end{align*}
	Thus, $\|T_2\|_{\mathscr H_2} \leq 1$, i.e. $T_2$ is similar to a contraction on $H_2.$
	An analogous argument shows that $T_1$ is similar to a contraction on $H_1$, so the ``only if'' part of the claim is proven for $m=2$.
	
	Next, employing the associative law for finite tensor products of Hilbert spaces and arguing as at the end of the proof of Theorem~\ref{split_finite}(ii), one shows that if the statement is true for all $k \leq m$, then it is also true for $k=m+1$, and the statement follows by induction on $m$. 
\end{proof}

To illustrate Theorem \ref{Discrete1}, as in the situation of one-parameter semigroups, we employ it 
to construct operators with specific properties that are nevertheless not similar to contractions.
Let $H$ be a Hilbert space and $T \in \linearOp(H)$. It was asked by  C. Davis whether the existence of the limits $\lim_{n\to \infty} \|T^n h\|$, $\lim_{n\to\infty} \|(T^\ast)^n h\|$ for every $h\in H$ implies that $T$ is similar to a contraction. This question was answered in the negative in \cite{eckstein1972contre-exemple} where, using  a version of Foguel's operator \cite{foguel1964counterexample}, it was constructed $T\in \linearOp(H)$ such that $T$  is not similar to a contraction and
$$\lim_{n\to \infty} \|T^n h\| = 0, \qquad \lim_{n\to\infty} \|(T^\ast)^n h\| =0, \qquad h \in H.
$$
Here we employ  Theorem~\ref{Discrete1} to give a different example of such an operator $T.$ 

\begin{proposition}
	There exists 
	$T\in \linearOp(H)$ such that $T$ is not similar to a contraction on $H$ and
	$$\lim_{n\to \infty} \|T^n h\| = 0, \qquad \lim_{n\to\infty} \|(T^\ast)^n h\| =0, \qquad h \in H.
	$$
\end{proposition}
\begin{proof}
	Let $F \in \linearOp(\ell^2(\NN) \oplus \ell^2(\NN))$ be any power bounded operator not similar to a contraction, e.g. Foguel's operator. Let $S$ be the unilateral shift on $\ell^2(\NN),$ and set $H:= \left(\ell^2(\NN) \oplus \ell^2(\NN)\right)  \otimes \ell^2(\NN)  \otimes \ell^2(\NN),$ and $T := F  \otimes S  \otimes S^\ast 
	\in \linearOp (H).$
	
	As $r(S) = r(S^\ast) =1$, it follows by Theorem~\ref{Discrete1} that $T$ is not similar to a contraction. On the other hand, as $\lim_{n\to\infty} (S^\ast)^n = 0$ strongly in $\linearOp(\ell^2(\mathbb N))$, it is readily seen that $\lim_{n\to\infty} T^n = 0$ strongly in $\mathcal L(H)$ and that $\lim_{n\to\infty} (T^\ast)^n = \lim_{n\to\infty} (F^\ast \otimes S^\ast  \otimes S)^n =0$ strongly in $\mathcal L(H),$ which finishes the proof. 
\end{proof}

\section{Similarities of tensor products and complete boundedness}

In this section we give an alternative proof of Theorem~\ref{split_finite}(ii), under the assumption that none of the factors is quasi-nilpotent. The proof relies on Paulsen’s similarity criterion, in the spirit of \cite[Theorem 2.3]{paulsen1995centered}.
In addition to Theorem~\ref{PPPTh}, we also require Propositions~\ref{contractiveStableProp} and~\ref{CrSimSetProp} to complete our proof.

First of all, we need to provide a few definitions concerning Paulsen's theorem. Let $H$ be a Hilbert space,
and let $\mathbb D$ stand for the unit disc. An operator $T$ in $\linearOp(H)$ is said to be \textit{polynomially bounded} if there exists $M\geq1$ such that
$$
\|p(T)\| \leq M \sup_{z\in \DD} |p(z)|, \qquad \mbox{for every polynomial } p.
$$
Similarly, $T$ is said to be \textit{completely polynomially bounded} if there exists $M \geq 1$ such that, for every $m \in \NN$ and every $m\times m$ matrix $P = \{p_{ij}\}_{i,j=1}^m$ of polynomials $p_{ij}$, one has
\begin{equation}\label{CPBeq}
	\|P(T)\| \leq M \sup_{z\in \DD} \|P(z)\|,
\end{equation}
  where $P(T)$ is the operator in $\bigoplus_{j=1}^m H$ given by $P(T) = \{p_{ij}(T)\}_{i,j=1}^m$. 
	If $T$ is completely polynomially bounded, then the smallest number $M \geq 1$ satisfying \eqref{CPBeq} is called the \textit{complete polynomial bound} of $T$, and we denote it by $M_{cpb}(T)$. If $T$ is not completely polynomially bounded, we set $M_{cpb}(T) = \infty$.

The next result due to Paulsen is classical and can be found e.g. in \cite{paulsen1984every}.
\begin{theorem}[Paulsen's theorem] \label{paul_th}
	Let $H$ be a Hilbert space and let $T \in \linearOp(H)$. Then $T$ is similar to a contraction if and only if $T$ is completely polynomially bounded. In addition, if $T$ is similar to a contraction, then $\mathcal C(T) = M_{cpb}(T)$.
\end{theorem}

Using Paulsen's theorem we obtain the following result, which was given in \cite[Theorem 2.3]{paulsen1995centered} in a slightly different context. Our proof proceeds along similar lines.  
\begin{theorem}\label{PPPTh}
	Let $H_1$ and $H_2$ be Hilbert spaces, and let $T_1 \in \linearOp(H_1)$ and $T_2 \in \linearOp(H_2)$ be such that $r(T_2) =1$ and $T_1 \otimes T_2$ is similar to a contraction on $H_1  \otimes H_2$. Then $T_1$ is similar to a contraction on $H_1$ and $\mathcal C(T_1) \leq \mathcal C(T_1  \otimes T_2)$.
\end{theorem}
\begin{proof}
	We may assume $1 \in \sigma(T_2)$ by replacing $T_2$ by $e^{i\theta}T_2$ for suitable $\theta \in [0,2\pi)$. Then $1$ belongs to the boundary of $\sigma(T_2)$ and is thus 
		an approximate eigenvalue of $T_2$. 
	Let $(f_n)_{n\in \NN} \subset H_2$ be such that $\|f_n\|=1$ for $n\in \NN$ and $\lim_{n\to \infty} \|T_2f_n - f_n\|=0$,
	hence for every $k \in \mathbb N,$
	\begin{equation*}
		\lim_{n\to \infty} \|T_2^k f_n - f_n\| = 0.
	\end{equation*}
	Fix $m \in \NN$ and let a family of polynomials $P=\{p_{ij}\}_{i,j=1}^m$ 
	be given by $p_{ij}(z) = \sum_{l=1}^\infty a_{ij}^l z^l$, $z \in \DD$, where $a_{ij}^l=0$ for all but finitely many $a_{ij}^l$. Furthermore, fix $\varepsilon > 0$  and $(h_1,  \ldots,h_m) \in \bigoplus_{j=1}^m H_1$ satisfying $\|P(T_1)(h_1, \ldots, h_m)\| \geq \|P(T_1)\| - \varepsilon$ and $\|(h_1,\ldots,h_m)\|=1$. Then $\|(h_1 \otimes f_n, \ldots, h_m \otimes f_n)\| = 1$, and for all $n \in \mathbb N,$
	\begin{align*}
		M_{cpb} (T_1 \otimes T_2)\sup_{z\in \DD} \|P(z)\| &\geq \|P (T_1\otimes T_2)\| \geq \|P (T_1\otimes T_2) (h_1 \otimes f_n, \ldots, h_m \otimes f_n)\|
		\\&= \left\|\left(\sum_{j=1}^m\sum_{l=1}^\infty \left(a_{1j}^l T_1^l h_j \otimes T_2^l f_n\right), 
		\ldots, \sum_{j=1}^m\sum_{l=1}^\infty \left( a_{mj}^l T_1^l h_j \otimes T_2^l f_n\right) \right) \right\|.
	\end{align*}  
	Therefore, \allowdisplaybreaks{
	\begin{align*}
		& M_{cpb} (T_1 \otimes T_2) \sup_{z\in \DD} \|P(z)\| 
		\\ \geq& \lim_{n\to \infty} \left\|\left(\sum_{j=1}^m\sum_{l=1}^\infty \left(a_{1j}^l T_1^l h_j \otimes T_2^l f_n\right),  \ldots, \sum_{j=1}^m\sum_{l=1}^\infty \left(a_{mj}^l T_1^l h_j \otimes T_2^l f_n \right) \right) \right\|
		\\ = & \lim_{n\to \infty} \left\|\left(\sum_{j=1}^m\sum_{l=1}^\infty \left(a_{1j}^l T_1^l h_j \otimes f_n\right),  \ldots, \sum_{j=1}^m\sum_{l=1}^\infty \left(a_{mj}^l T_1^l h_j \otimes f_n \right)\right) \right\|
		\\ = & \lim_{n\to \infty} \left\|\left(\left(\sum_{j=1}^m\sum_{l=1}^\infty a_{1j}^l T_1^l h_j\right) \otimes f_n,  \ldots, \left(\sum_{j=1}^m\sum_{l=1}^\infty a_{mj}^l T_1^l h_j\right) \otimes f_n \right)  \right\|
		\\ = & \left\|\left( \sum_{j=1}^m\sum_{l=1}^\infty a_{1j}^l T_1^l h_j,  \ldots, \sum_{j=1}^m\sum_{l=1}^\infty a_{mj}^l T_1^l h_j \right) \right\|
		\\ = & \|P(T_1) (h_1, \ldots, h_m)\| \geq \|P(T_1)\|-\varepsilon.
	\end{align*}
	As $\varepsilon>0$, $m \in \NN$ and $P$ were arbitrary, we conclude that $T_1$ is completely polynomially bounded and $M_{cpb} (T_1) \leq M_{cpb} (T_1 \otimes T_2)$, so Theorem \ref{paul_th} implies the claim.} 
\end{proof}

Theorem~\ref{Discrete1} 
 is then a direct consequence 
 of Theorem~\ref{PPPTh}, since the general case $m \in \NN$ follows by induction on $m$ once we have proven the statement for $m=2$ 
as in the proof of Theorem~\ref{PPPTh}.

We also need the following result proved in \cite[Proposition 6.1]{oliva2025similarity}
and inspired by \cite[Proposition 2.5]{holbrook1977distortion}. 
Since \cite{oliva2025similarity} is very recent, and to be self-contained we provide
a short argument.

\begin{proposition}\label{CrSimSetProp}
	Let $H$ be a Hilbert space, and let $\mathcal T = (T(t))_{t\geq0} \subset \linearOp(H)$ be a semigroup strongly continuous in $(0,\infty)$. Then $\mathcal T \in \mathcal{SC}(H)$ if and only if $\liminf_{t\to0} \mathcal C(T(t)) < \infty$. If the latter condition holds, then
	$$\mathcal C(\mathcal T) = \lim_{t\to 0} \mathcal C(T(t)) = \sup_{t>0} \mathcal C(T(t)).
	$$
\end{proposition}
\begin{proof}
	The ``only if'' part of the statement is trivial since $\mathcal C(T(t)) \leq \mathcal C(\mathcal T)$ for all $t\geq0$. Thus, we assume that $\mathcal C :=\liminf_{t\to0} \mathcal C(T(t)) < \infty$ and fix $(t_n)_{n\in \NN} \subset (0,\infty)$ such that $\lim_{n\to \infty} t_n = 0$ and 
	$$
	\lim_{n \to \infty} \mathcal C(T(t_n)) = \mathcal C.
	$$
	For each $n \in \NN$, let $\|\cdot\|_n$ be a Hilbertian norm on $H$ satisfying $\|T(t_n)\|_{n} \leq 1$ and
	\begin{align*}
		\|h\| \leq \|h\|_n \leq \mathcal C(T(t_n)) \|h\|, \qquad h\in H.
	\end{align*}
	Define 
	$$
	\|h\|_{\rm{eq}}^2 := \operatorname{LIM}   \bigl[\|h\|_n^2 \bigr], \qquad h\in H,
	$$
	and note that $\|\cdot\|_{\rm{eq}}$ is a Hilbertian norm on $H$ such that
	$$\|h\| \leq \| h \|_{\rm{eq}} \leq  \mathcal C \|h\|, \qquad h \in H.
	$$
	Now fix $t>0$ and $h \in H$. For each $n\in \NN$, let $m_n \in \NN$ be such that $m_n t_n$ realizes the distance between $\{t\}$ and $t_n \NN$. Then $\lim_{n\to \infty} |t-m_n t_n| \leq \lim_{n\to \infty} t_n = 0$ and
	$$
	\|T(t)h - T(m_n t_n)h\|_n \leq \mathcal C(T(t_n)) \|T(t) h - T(m_n t_n)h\| 
	\rightarrow 0, \qquad \text{as}\,\, n \to \infty.
	$$
	In addition, $\|T(m_n t_n)\|_n \leq \|T(t_n)\|_n^{m_n} \leq 1$, $n \in \NN$. Therefore,
	\begin{align*}
		\|T(t)h\|_{\rm{eq}}^2 &= \operatorname{LIM} \bigl[\|T(t)h\|_n^2 \bigr] \leq 
		\operatorname{LIM}  \left[\left(\|T(m_n t_n)h\|_n + \|T(t)h - T(m_n t_n)h\|_n \right)^2\right]
		\\ &= \operatorname{LIM}  \bigl[\|T(m_n t_n) h\|_n^2 \bigr] 
		\leq \operatorname{LIM}   \bigl[\|h\|_n^2 \bigr] = \|h\|_{\rm{eq}}^2, \qquad t \geq 0, \, h \in H.
	\end{align*}
	Thus,  $\|T(t)\|_{\rm{eq}} \leq 1$ for all $t\geq0,$ so that 
	$\mathcal T \in \mathcal{SC}(H)$ and $\mathcal C(\mathcal T) \leq \mathcal C$. 
	Since 
	$$
	\mathcal C(\mathcal T) \geq \sup_{t>0} \mathcal C(T(t)) \geq \liminf_{t\to0} 
	\mathcal C(T(t)) = \mathcal C,
	$$
	we have
	$$
	\mathcal C(\mathcal T) = \sup_{t>0} \mathcal C(T(t)) = \liminf_{t\to0} \mathcal C(T(t)).
	$$
	Hence $\lim_{t\to0} \mathcal C(T(t)) = \mathcal C(\mathcal T)$, and the proof is finished.
\end{proof}
\medskip

Now we are able to provide an alternative proof of
Theorem~\ref{split_finite}(ii).

\begin{proof}[Proof of Theorem \ref{split_finite}(ii) via complete boundedness]
The ``if'' part of the statement can be proved by a direct argument similarly 
to the ``if'' part of the proof of  Theorem \ref{split_finite}(ii),
so we concentrate on the ``only if'' part.

	Arguing as at the end of the proof of Theorem~\ref{split_finite}(ii), the statement for all $m\in \NN$ follows by induction once we have proven it for $m=2$. To consider this case, let $H_1$ and $H_2$ be Hilbert spaces and let $\mathcal T_1 = (T_1(t))_{t\geq0} \subset \linearOp(H)$ and 
	$\mathcal T_2 = (T_2(t))_{t\geq 0} \subset \linearOp(H_2)$ be semigroups that are not quasi-nilpotent and satisfying 
	the assumptions of Theorem \ref{split_finite}(ii). Set $a = -\omega_0(\mathcal T_2) \in \RR$, so that $\omega_0(e_{a} \mathcal T_2) = \omega_0(\mathcal T_2) +a =0$ and $r(e^{a t} T_2(t))=1$ for every $t\geq0$. 
	By Theorem~\ref{PPPTh}, $e^{-at}T_1(t)$ is similar to a contraction for every $t\geq0$ and, moreover,
	$$\mathcal C(e^{-a t}T_1(t)) \leq \mathcal C(e^{-a t} T_1(t)\otimes e^{a t} T_2(t))   \leq \mathcal C(\mathcal T_1 \otimes \mathcal T_2), \qquad t\geq0.
	$$
	Therefore, it follows from Proposition~\ref{CrSimSetProp} that $e_{-a} \mathcal T_1 \in \mathcal{SC}(H_1)$. In particular, $\mathcal T_1 \in \mathcal{SQC}(H_1)$. An analogous argument shows that 
	$e_{-b} \mathcal T_2  \in \mathcal{SC}(H_2)$ where $b = -\omega_0(\mathcal T_1)$. 
	
	If $\omega_0(\mathcal T_1 \otimes \mathcal T_2)= \omega_0(\mathcal T_1) + \omega_0(\mathcal T_2)  =0$, then $b = -a$ and the statement follows by what we have already shown. If $\omega_0(\mathcal T_1 \otimes \mathcal T_2) < 0$, then choose $d \in \left(\omega_0(\mathcal T_2), -\omega_0(\mathcal T_1)\right)$ so that $\omega_0(e_d \mathcal T_1) = \omega_0(\mathcal T_1) + d < 0$ and $\omega_0(e_{-d} \mathcal T_2) = \omega_0(\mathcal T_2) - d < 0$. 
	Thus, $e_d \mathcal  T_1 \in \mathcal{SQC}(H_1)$ and $e_{-d} \mathcal T_2 \in \mathcal{SQC}(H_2)$ and both  $e_d \mathcal  T_1$ and $e_{-d} \mathcal T_2$ are exponentially stable. Hence $e_d \mathcal  T_1 \in \mathcal{SC}(H_1)$ and $e_{-d} \mathcal T_2 \in \mathcal{SC}(H_2)$ by Proposition~\ref{contractiveStableProp}, and the proof is finished.
\end{proof}

\begin{remark} 
	 Using the same notation as in Theorem~\ref{split_finite}(ii), we could not address the case where $\omega_0(\mathcal T_k)=-\infty$ (i.e. $\mathcal T_k$ is a quasi-nilpotent semigroup) for some $k\in \mathbb N_m$, using Paulsen's theorem and Theorem~\ref{PPPTh}.
		Recall that this case, treated in Theorem~\ref{split_finite}(ii) by different techniques, is crucial in Section~\ref{counterexamplesSub} for constructing quasi-nilpotent semigroups (among other properties) not similar to contractions.
\end{remark}

\bibliography{mybib}
\bibliographystyle{amsplain-nodash}  
\end{document}